\documentclass[12pt]{amsart}
\usepackage{amssymb}
\usepackage{amsmath}
\usepackage{longtable}
\newcommand\g{{\mathfrak g}}
\newcommand\h{{\mathfrak h}}

\renewcommand{\v}{\mathfrak{v}}

\newcommand\gl{\mathfrak{gl}}
\renewcommand{\a}{\mathfrak{a}}

\newcommand\m{\mathfrak m}
\newcommand\lfr{\mathfrak l}
\newcommand\CS{\mathcal{CS}}
\newcommand\n{\mathfrak n}
\newcommand\q{\mathfrak q}
\newcommand\z{\mathfrak z}
\newcommand\s{\mathfrak s}

\renewcommand{\t}{\mathfrak{t}}

\newcommand\im{\operatorname{im}}

\newcommand\Pic{\operatorname{Pic}}

\newcommand\Aut{\operatorname{Aut}}

\newcommand\Hom{\operatorname{Hom}}

\newcommand\Gr{\operatorname{Gr}}
\newcommand\Supp{\operatorname{Supp}}
\newcommand\K{{\mathbb{K}}}
\newcommand\X{\mathfrak X}
\newcommand\Q{\mathbb Q}

\newcommand\D{\mathcal D}
\newcommand\Z{\mathbb Z}
\newcommand\A{\mathfrak A}
\newcommand\V{\mathcal V}

\renewcommand\H{\mathcal{H}}
\newcommand\Rad{\operatorname{R}}

\newcommand\so{\mathfrak{so}}
\renewcommand\sp{\mathfrak{sp}}
\newcommand\spin{\mathfrak{spin}}
\newcommand\ord{\operatorname{ord}}
\newcommand\GL{\mathop{\rm GL}\nolimits}

\newcommand\SL{\mathop{\rm SL}\nolimits}
\newcommand\Sp{\mathop{\rm Sp}\nolimits}
\newcommand\SO{\mathop{\rm SO}\nolimits}

\newcommand\Span{\operatorname{Span}}

\newcommand{\Ad}{\mathop{\rm Ad}\nolimits}

\newcommand{\rank}{\mathop{\rm rk}\nolimits}

\renewcommand{\Ad}{\mathop{\rm Ad}\nolimits}
\newcommand\quo{/\!/}
\newtheorem{Thm}{Theorem}
\newtheorem{Prop}{Proposition}[subsection]
\newtheorem{Cor}[Prop]{Corollary}

\newtheorem{Lem}[Prop]{Lemma}
\theoremstyle{definition}
\newtheorem{Ex}[Prop]{Example}
\newtheorem{defi}[Prop]{Definition}
\newtheorem{Rem}[Prop]{Remark}
\textheight=232mm \numberwithin{equation}{section}
\author{Ivan V. Losev}
\title{Uniqueness property for spherical homogeneous spaces}
\thanks{Partially supported by A. Moebius foundation}
\begin{document}
\begin{abstract}
Let $G$ be a connected reductive group. Recall that a homogeneous
$G$-space $X$ is called spherical if  a Borel subgroup $B\subset G$
has an open orbit on $X$. To $X$ one assigns certain combinatorial
invariants: the weight lattice, the valuation cone and the set of
$B$-stable prime divisors. We prove that two spherical homogeneous
spaces with the same combinatorial invariants are equivariantly
isomorphic. Further, we recover the group of $G$-equivariant
automorphisms of $X$ from these invariants.
\end{abstract}
\maketitle
\section{Introduction}
At first, let us fix some notation and introduce some terminology.

Throughout the paper the base field $\K$ is  algebraically closed
and of characteristic zero.   Let $G$ be a connected reductive
group. Fix a Borel subgroup $B\subset G$. An irreducible $G$-variety
$X$ is said to be {\it spherical} if $X$ is normal and $B$ has an
open orbit on $X$. The last condition is equivalent to
$\K(X)^{B}=\K$. Note that a spherical $G$-variety contains an open
$G$-orbit. An algebraic subgroup $H\subset G$ is called spherical if
 $G/H$ is spherical.

The theory of spherical varieties and related developments seem to
be the most important topic  in the study of  algebraic
transformation groups for the last twenty five years (see
\cite{Brion_rev}, \cite{Knop5}, \cite{Timashev_rev} for a review of
this theory). The first problem arising in the study of spherical
varieties is their classification. The basic result here is the
Luna-Vust theory, \cite{LV}, describing all spherical $G$-varieties
with a given open $G$-orbit. The description is carried out in terms
of certain combinatorial invariants. So the classification of all
spherical varieties is reduced to the description of all spherical
subgroups and the computation of the corresponding combinatorial
invariants. However, there are  many spherical subgroups and
their explicit classification (i.e.,  as subsets of $G$)
 seems to be possible only in some very special cases.
For example, connected reductive spherical subgroups were classified
(in fact, partially) in \cite{Kramer}, \cite{Brion},\cite{Mikityuk}.
In the general case one should perform the classification using a
different language. A reasonable language was proposed by Luna,
\cite{Luna5}. Essentially, his idea is to describe spherical
homogeneous spaces in terms of combinatorial invariants established
in \cite{LV}. One should prove that these invariants determine a
spherical homogeneous space uniquely  and then check the existence
of the  spherical homogenous space corresponding to any set of
invariants satisfying some combinatorial conditions. Luna completed
his program in \cite{Luna5} for groups of type $A$ (that is, when
any simple normal subgroup is locally isomorphic to a special linear
group). Later on, the full classification for groups of type $A-D$
was carried out by Bravi, \cite{Bravi}. A partial classification for
type $A-C$ was obtained in \cite{Pezzini}. Finally, recently Bravi
obtained the classification in types A-D-E, \cite{Bravi2}. The main
result of the present paper is the proof of the uniqueness part of
Luna's program. An advantage of our approach over Luna's is that we do not use
long case-by-case considerations.

Let us describe combinatorial invariants in interest. We will use standard
notation  recalled in Section 2 below. Fix a maximal
torus $T\subset B$.

Let $X$ denote a spherical $G$-variety.  The  set
$\X_{G,X}:=\{\mu\in\X(T)| \K(X)^{(B)}_\mu\neq \{0\}\}$ is called the
{\it weight lattice} of $X$. This is a sublattice in $\X(T)$. By the
{\it Cartan space} of $X$ we mean $\a_{G,X}:=\X_{G,X}\otimes_\Z\Q$.
This is a subspace in $\t(\Q)^*$.

Next we define the valuation cone of $X$. Let $v$ be a $\Q$-valued
discrete $G$-invariant valuation of $\K(X)$. Since $X$ is
 spherical, we have
$\dim\K(X)^{(B)}_\mu=1$ for any $\mu\in \X_{G,X}$. So one defines the element
$\varphi_v\in \a_{G,X}^*$ by the formula
\begin{equation*}
\langle\varphi_v,\mu\rangle=v(f_\mu), \forall \mu\in\X_{G,X},
f_{\mu}\in \K(X)^{(B)}_\mu\setminus \{0\}.
\end{equation*}

  It is known, see \cite{Knop5}, that the map $v\mapsto \varphi_v$
 is injective. Its image is a finitely generated
convex cone in $\a_{G,X}^*$.
 We denote this cone by $\V_{G,X}$
and call it the {\it valuation cone} of $X$.

Let $\D_{G,X}$  denote the set of all prime $B$-stable divisors of
$X$. This is a finite set. To $D\in \D_{G,X}$ we assign
$\varphi_D\in \a_{G,X}^*$ by $\langle\varphi_D,\mu\rangle=\ord_D
(f_\mu), \mu\in \X_{G,X}, f_\mu\in \K(X)^{(B)}_\mu\setminus\{0\}$.
Further, for $D\in D_{G,X}$ set $G_D:=\{g\in G| gD=D\}$. Clearly,
$G_D$ is a parabolic subgroup of $G$ containing $B$. Choose
$\alpha\in \Pi(\g)$. Set $\D_{G,X}(\alpha):=\{D\in G_D|
P_\alpha\not\subset G_D\}$. Here and below by $P_\alpha$ we denote
the parabolic subgroup of $G$ generated by $B$ and the
one-dimensional unipotent subgroup of $G$ corresponding to the root
subspace of weight $-\alpha$.

Below we regard $\D_{G,X}$ as an abstract set equipped with two maps
$D\mapsto \varphi_D, D\mapsto G_D$. For instance, if $X_1,X_2$ are
spherical $G$-varieties, then, when writing $\D_{G,X_1}=\D_{G,X_2}$,
we mean that there exists a bijection $\iota:\D_{G,X_1}\rightarrow
\D_{G,X_2}$ such that $G_D=G_{\iota(D)},
\varphi_D=\varphi_{\iota(D)}$.

\begin{Thm}\label{Thm:1.1}
Let $H_1,H_2$ be spherical subgroups of $G$. If
$\X_{G,G/H_1}=\X_{G,G/H_2},\V_{G,G/H_1}=\V_{G,G/H_2},\D_{G,G/H_1}=\D_{G,G/H_2}$,
then $H_1,H_2$ are $G$-conjugate.
\end{Thm}

Thanks to Theorem \ref{Thm:1.1}, one may hope to describe all
invariants of a spherical homogeneous space $X$ in terms of
$\X_{G,X},\V_{G,X},\D_{G,X}$.

We want to describe the group $\Aut^G(X)$ of $G$-equivariant
automorphisms of $X$.   For
any $\varphi\in \Aut^G(X), \lambda\in \X_{G,X}$ there is $a_{\varphi,\lambda}\in
\K^\times$ such that
$\varphi|_{\K(X)^{(B)}_\lambda}=a_{\varphi,\lambda}id$. The map
$(\lambda,\varphi)\mapsto a_{\varphi,\lambda}$ gives rise to the
homomorphism $\Aut^G(X)\rightarrow
A_{G,X}:=\Hom_\Z(\X_{G,X},\K^\times)$. This homomorphism is
injective and its image  $\A_{G,X}$ is closed, see, for example,
\cite{Knop8}, Theorem 5.5. Therefore $\A_{G,X}$ is recovered from
$\Lambda_{G,X}:=\{\lambda\in\X_{G,X}|
\langle\lambda,\A_{G,X}\rangle=1\}$. We call $\Lambda_{G,X}$ the
{\it root lattice} of $X$. By the root lattice of an arbitrary
spherical $G$-variety $X$ (also denoted by $\Lambda_{G,X}$) we mean
the root lattice of the open $G$-orbit in $X$.

To describe $\Lambda_{G,X}$ in terms of
$\X_{G,X},\mathcal{V}_{G,X},\D_{G,X}$ we need to recall some further
facts about the structure of $\mathcal{V}_{G,X}$. Namely, fix an
$N_G(T)$-invariant scalar product on  $\t(\Q)$. It induces the
scalar product on $\a_{G,X}$. With respect to this scalar product
$\V_{G,X}$ becomes a Weyl chamber for a (uniquely determined) linear
group $W_{G,X}$ generated by reflections, see \cite{Brion}.  It
turns out, see, for example, \cite{Knop8}, Sections 4,6, that
$\Lambda_{G,X},\X_{G,X}$ are $W_{G,X}$-stable. Denote by
$\Psi_{G,X}$ (resp., $\overline{\Psi}_{G,X}$) the set of primitive
elements $\alpha\in \X_{G,X}$ (resp., $\alpha\in \Lambda_{G,X}$)
such that $\ker\alpha\subset \a_{G,X}^*$ is a wall of $\V_{G,X}$ and
$\langle\alpha,\V_{G,X}\rangle\leqslant 0$. Clearly,
$\Psi_{G,X},\overline{\Psi}_{G,X}$ are systems of simple roots with
Weyl group $W_{G,X}$. An element of $\Psi_{G,X}$ is called a {\it
spherical root} of $X$. It was proved in \cite{Knop8} that
$\overline{\Psi}_{G,X}$ generates $\Lambda_{G,X}$. On the other
hand, in general,
$\Psi_{G,X}$  does not  generate $\X_{G,X}$. 

We recover $\overline{\Psi}_{G,X}$ from $\Psi_{G,X},\D_{G,X}$. To
do this  we define a certain  subset $\Psi_{G,X}^+\subset
\Psi_{G,X}$ of {\it distinguished} roots, see Definition
\ref{defi:1}.

\begin{Thm}\label{Thm:1.2}
$\overline{\Psi}_{G,X}= (\Psi_{G,X}\cap \Lambda(\g)\setminus
\Psi_{G,X}^+)\sqcup\{2\alpha|\alpha\in \Psi_{G,X}^+\cup
(\Psi_{G,X}\setminus \Lambda(\g))\}$.
\end{Thm}

There are other invariants of a spherical homogeneous space $G/H$
that can be described in terms of
$\X_{G,G/H},\V_{G,G/H},\D_{G,G/H}$. For example, Knop, \cite{Knop5},
described  the set $\H_H$ of all algebraic subgroups
$\widetilde{H}\subset G$ such that $H\subset \widetilde{H}$ and
$\widetilde{H}/H$ is connected, see Subsection
\ref{SUBSECTION_inclusions} for details. This description plays a
crucial role in our proofs.

Theorems \ref{Thm:1.1},\ref{Thm:1.2} seem to be quite independent
from each other. The reason why they are brought together in a
single paper is twofold. First, the key ideas of their proofs are
much alike. Second, we use Theorem \ref{Thm:1.2} in the proof of
Theorem \ref{Thm:1.1}.

Now let us briefly describe the content of this paper. In Section
\ref{SECTION_Notation} we present conventions and the list of
notation we use. In Section \ref{SECTION_prelim} we gather some
preliminary results. Most of them are standard. Section
\ref{SECTION_Proof2} is devoted to the proofs of Theorems
\ref{Thm:1.1},\ref{Thm:1.2}.  In the beginning of Sections
\ref{SECTION_prelim},\ref{SECTION_Proof2} their content is described
in more detail. All propositions, definitions, etc. are numbered
within a subsection.

Finally, we would like to discuss applications of our results.
Theorem \ref{Thm:1.1}  significantly simplifies the proof of the
Knop conjecture, \cite{Knop_conj},  allowing  to get rid of many
ugly technical considerations.

As for Theorem \ref{Thm:1.2}, we   use it in \cite{Losev_Dem} to
prove  Brion's conjecture on the smoothness of Demazure embeddings,
see \cite{Brion}. In fact, it is that application that motivated us
to state Theorem \ref{Thm:1.2}. Let us recall the definition of the
Demazure embedding. Let $\h$ be a subalgebra of $\g$. Suppose
$\h=\n_\g(\h)$ and  the subgroup $H:=N_G(\h)\subset G$ is spherical.
Consider $\h$ as a point of the Grassmanian $\Gr_d(\g)$. The
$G$-orbit of $\h$ is isomorphic to the homogeneous space $G/H$. By
the Demazure embedding of $G/H$ we mean the closure $\overline{G\h}$
of $G\h$ in $\Gr_d(\g)$. Brion conjectured that $\overline{G\h}$ is
smooth. Again, for groups of type $A$ this conjecture was proved by
Luna in \cite{Luna_Dem}.

{\bf Acknowledgements} This paper was written during author's stay
in Rutgers University, New Brunswick,  in the beginning of 2007. I
am grateful to this institution and especially to Professor F. Knop
for hospitality. I also would like to thank  F. Knop for stimulating
discussions and inspiration. Finally, I thank the referees for useful remarks
on a previous version of this paper.

\section{Notation, terminology and
conventions}\label{SECTION_Notation}  If an algebraic group is
denoted by a capital Latin letter, then we denote its Lie algebra by
the corresponding small German letter.

We fix a Borel subgroup $B\subset G$ and a maximal torus $T\subset
B$. This allows us  to define the root system $\Delta(\g)$, the Weyl
group $W(\g)$, the root lattice $\Lambda(\g)$ and the system of
simple roots $\Pi(\g)$ of $\g$.  For simple roots and fundamental
weights  we use the notation of \cite{VO}. Note that the character
groups $\X(T),\X(B)$ are naturally identified. By $B^-$ we denote
the Borel subgroup of $G$ opposite to $B$ and containing $T$.

By a spherical subalgebra of $\g$ we mean the Lie algebra of a
spherical subgroup of $G$.


If $X_1,X_2$ are  $G$-varieties, then we write  $X_1\cong^G X_2$
when $X_1,X_2$ are $G$-equivariantly isomorphic.  Note that if
$V_1,V_2$ are $G$-modules, then $V_1\cong^G V_2$ iff $V_1,V_2$ are
isomorphic as $G$-modules.

Let $X_1,X_2$ be spherical varieties.  We write $X_1\equiv^G X_2$ if
$\X_{G,X_1}=\X_{G,X_2}, \V_{G,X_1}=\V_{G,X_2},
\D_{G,X_1}=\D_{G,X_2}$. When $X_1=G/H_1,X_2=G/H_2$ and $X_1\equiv^G
X_2$ we write $H_1\equiv^G H_2$.

Let $Q$ be a parabolic subgroup of $G$ containing either $B$ or
$B^-$. There is a unique Levi subgroup of $Q$ containing $T$, we
call it the {\it standard} Levi subgroup of $Q$.

Now let $S$ be a reductive group and $\Gamma$ be a  group with a
fixed homomorphism to the group of outer automorphisms of $S$. There
is a natural right action of $\Gamma$ on the set of (isomorphism
classes of) $S$-modules. Namely, let $V$ be an $S$-module and
$\rho:G\rightarrow \GL(V)$ be the corresponding representation. For
$V^\gamma, \gamma\in \Gamma,$ we take  the $S$-module corresponding
to the representation $\rho\circ\overline{\gamma}$, where
$\overline{\gamma}$ is a representative of the image of $\gamma$.

\begin{longtable}{p{2cm}p{14cm}}
$\sim_\Gamma$& the equivalence relation induced by an action of a
group $\Gamma$\\
$A^{(B)}_\lambda$&$=\{a\in A| b.a=\lambda(b)a, \forall b\in B\}$.
\\$A^{(B)}$&$=\cup_{\lambda\in \X(B)}A^{(B)}_\lambda$.
\\ $A_{G,X}$&$=\Hom(\X_{G,X},\K^\times)$.\\
$\a_{G,X}$& the Cartan space of a spherical $G$-variety $X$.
\\ $\A_{G,X}$& the image of $\Aut^G(X^0)$ in $A_{G,X}$, where $X^0$
is the open $G$-orbit of a spherical variety $X$.\\
$\A_{G,X}(\alpha)$&$=\{\varphi\in \A_{G,X}|
\langle\alpha,\varphi\rangle=-1\}$.
\\ $(\a_H^{\widetilde{H}},\D_H^{\widetilde{H}})$& the colored
subspace  associated with $\widetilde{H}\in \H_H$.
\\ $\Aut^G(X)$& the group of all
$G$-equivariant automorphisms of  $X$.
\\ $\CS_{G,X}$& the set of all colored subspaces of $(\a_{G,X}^*,\D_{G,X})$
\\ $\D_{G,X}$& the set of all prime $B$-stable divisors of  $X$.
\\ $\D_{G,X}(\alpha)$& $=\{D\in\D_{G,X}| P_\alpha\not\subset G_D\}$.
\\ $f_\lambda$& a nonzero element in $\K(X)^{(B)}_\lambda$.
\\ $(f)$& the  divisor of a rational function $f$.
\\ $(G,G)$& the derived subgroup of a group
$G$\\
$[\g,\g]$& the derived subalgebra of a Lie algebra $\g$.
\\ $G^{\circ}$& the unit component  of an algebraic
group $G$.
\\ $G*_HV$& the homogeneous bundle over $G/H$ with a fiber $V$.
\\ $G_x$& the stabilizer of $x$ under an action
of $G$.
\\ $\Gr_d(V)$& the Grassman variety consisting of all
$d$-dimensional subspaces of a vector space $V$.
\\ $\H_H$& the set of all algebraic subgroups $\widetilde{H}$ of  $G$ such that $H\subset \widetilde{H}$ and $\widetilde{H}/H$ is
connected.
\\ $\underline{\H}_H$&$=\{\widetilde{H}\in \H_H| \Rad_u(H)\subset \Rad_u(\widetilde{H}),
\widetilde{H}/\Rad_u(\widetilde{H})=H/\Rad_u(H),\newline
\Rad_u(\widetilde{\h})/\Rad_u(\h)\text{ is an irreducible }
H\text{-module}\}$.
\\ $\H^{\widetilde{H}}$& the set of all algebraic subgroups $H\subset \widetilde{H}$  such that
$\widetilde{H}/H$ is connected.
\\ $\overline{\H}^{\widetilde{H}}$&$=\{H\in \H^{\widetilde{H}}| \Rad_u(H)\subset
\Rad_u(\widetilde{H})\}$.
\\ $\K^\times$& the one-dimensional torus.
\\ $L_{G,X}$& the standard Levi subgroup of $P_{G,X}$.
\\ $N_G(H)$& $=\{g\in G| gHg^{-1}=H\}$\\ $N_G(\h)$& $=\{g\in G| \Ad(g)\h=\h\}$.\\
$\n_\g(\h)$& $=\{\xi\in\g| [\xi,\h]\subset\h\}$.
\\ $P_{G,X}$&$=\bigcap_{D\in \D_{G,X}}G_D$.
\\ $P_\alpha$& the minimal parabolic subgroup of $G$ containing
$B$ associated with $\alpha\in \Pi(\g)$.
\\ $\Pic(X)$& the Picard group of a variety $X$.
\\ $\Pic_G(X)$& the equivariant Picard group of a $G$-variety $X$.
\\ $\Rad_u(G)$& the unipotent radical of an algebraic group $G$.
\\ $\rank_G(X)$& $=\rank\X_{G,X}$.
\\ $\Span_{A}(M)$&$=\{a_1m_1+\ldots+a_km_k, a_i\in A, m_i\in M\}$.
\\ $\Supp(\gamma)$& the support of  $\gamma\in\Span_\Q(\Pi(\g))$, that is, the set $\{\alpha\in \Pi(\g)|n_\alpha\neq 0\}$, where
$\gamma=\sum_{\alpha\in\Pi(\g)}n_\alpha\alpha$.
\\
$\X(G)$& the character group of an algebraic group $G$.\\
$\X_{G,X}$& the weight lattice of a spherical $G$-variety $X$.
\\  $\X^+_{G,X}$&
the weight monoid of a spherical $G$-variety $X$.\\
$X/\Gamma$& the geometric quotient for an action $\Gamma:X$.\\
 $\#X$& the cardinality of a
set $X$.
\\ $\V_{G,X}$& the valuation cone of a spherical $G$-variety $X$.
\\$V(\mu)$& the irreducible module with  highest weight $\mu$.
\\ $W(\g)$& the Weyl group of a reductive Lie algebra $\g$.
\\   $Z_G(\h)$&=$\{g\in G| \Ad(g)|_\h=id\}$\\  $\z_\g(\h)$&$=\{\xi\in\g| [\xi,\h]=0\}$.
\\ $Z(G)$&the center of a group $G$.\\
$\z(\g)$&$=\z_\g(\g).$
\\ $\alpha^\vee$& the dual root corresponding to a root $\alpha$.\\
$\langle\alpha,v\rangle$& the pairing of elements $\alpha,v$ of dual
vector spaces or dual abelian groups.
\\  $\Delta(\g)$& the root system of  $\g$.
\\ $\Lambda(\g)$& the root lattice of  $\g$.
\\ $\Lambda_{G,X}$& the root lattice of $X$.
\\ $\Pi(\g)$& the system of simple roots of $\g$.
\\ $\Psi_{G,X}$& the system of spherical roots of a spherical
$G$-variety $X$.
\\ $\Psi_{G,X}^+$& the set  of all
distinguished elements of $\Psi_{G,X}$.
\\ $\Psi_{G,X}^i$& the subset of $\Psi_{G,X}^+$ of all
roots of type $i$, $i=1,2,3$.
\\ $\varphi_D$& the vector in $\a_{G,X}^*$ associated with $D\in
\D_{G,X}$.
\end{longtable}

\section{Preliminaries}\label{SECTION_prelim}
This section does not contain new results. Its goal is to recall
some basic facts about spherical varieties used in the proofs of the
main theorems. We often present short proofs of some results when we
do not know an appropriate reference.

 $X$ denotes a spherical $G$-variety.
 Throughout the section we set
$\X:=\X_{G,X},\D:=\D_{G,X},\Psi:=\Psi_{G,X},P:=P_{G,X}$, etc.

In Subsection \ref{SUBSECTION_auto} we recall some basic facts about
the group $\Aut^G(X)$. Subsection \ref{SUBSECTION_wonderful} is
devoted to spherical homogeneous spaces  admitting a so called
wonderful embedding. We also present there results concerning
systems of spherical roots. In Subsection
\ref{SUBSECTION_finiteness} we state a finiteness result for
spherical subalgebras coinciding with their normalizers. It plays an
important role in the proofs of Theorems
\ref{Thm:1.1},\ref{Thm:1.2}.

Subsection \ref{SUBSECTION_inclusions} deals with the structure of
the set $\H_H$. This description is due to Knop, \cite{Knop5}, and
is of major importance in our proofs. We describe the structure of
the ordered set  $\H_H$ in terms of so called {\it colored}
subspaces, see Definition \ref{defi:2}. Then we describe the
combinatorial invariants of $G/\widetilde{H},\widetilde{H}\in \H_H$.

Subsection \ref{SUBSECTION_incl_parab} is devoted to different
results related to the local structure theorem. This theorem is used
to reduce the study of the action $G:X$ to  the study of an action
of a certain Levi subgroup  $M\subset G$ on a certain subvariety
$X'\subset X$. Mostly, we use the theorem in the situation described
in Example \ref{Ex:2.5}. Originally, the local structure theorem was
proved independently in \cite{BLV},\cite{Grosshans}, but we use its
variant due to Knop, \cite{Knop3}.

Finally, in Subsection \ref{SUBSECTION_affine_prelim} we gather
different simple properties of affine spherical varieties.

\subsection{Equivariant automorphisms}\label{SUBSECTION_auto}
Throughout this subsection $X=G/H$.

One can identify  $\Aut^G(X)$  with $N_G(H)/H$. The group
$\Aut^G(X)$ acts naturally on $\D$ and the maps $D\mapsto G_D,
D\mapsto \varphi_D$ are $\Aut^G(X)$-invariant.

As it was explained in Introduction, there is a natural monomorphism
$\Aut^G(X)\rightarrow A$. So $\Aut^G(X)$ is commutative.

\begin{Lem}\label{Cor:2.4}
$N_G(H)=N_G(H^\circ)$.
\end{Lem}
\begin{proof}
Since $H^\circ$ is also spherical, the group $N_G(H^\circ)/H^\circ$
is commutative. Thus $H$ is a normal subgroup of $N_G(H^\circ)$.
\end{proof}

The next lemma seems to be well-known, see, for example,
\cite{Losev}, Lemma 7.17.

\begin{Lem}\label{Lem:2.2}
Let $\overline{X}$ be a quasiaffine spherical variety, whose open
$G$-orbit is isomorphic to $X$. Then
$\Aut^G(\overline{X})=\Aut^G(X)$.
\end{Lem}

\begin{Lem}\label{Lem:2.11}
Let $X_1,X_2$ be quasiaffine spherical $G$-varieties and
$\psi:X_1\rightarrow X_2$ a dominant $G$-equivariant morphism. For any
$\varphi\in \Aut^G(X_1)$ there exists a unique element
$\underline{\varphi}\in \Aut^G(X_2)$ such that
$\psi\circ\varphi=\underline{\varphi}\circ\psi$.
\end{Lem}
\begin{proof}
Note that $\psi^*(\K[X_2])\subset \K[X_1]$ is $\varphi$-stable. For
$\underline{\varphi}$ we take a unique element of $\Aut^G(X_2)$
coinciding with $\varphi$ on $\K[X_2]$.  \end{proof}

\begin{Lem}[\cite{Knop8}, Corollary 6.5]\label{Lem:2.1}
The Lie algebra of $\A$ coincides with $\K\otimes_{\Q}(\V\cap -\V)$.
\end{Lem}


\begin{Lem}\label{Lem:2.1.2}
Let $\Gamma$ be an algebraic subgroup of $\A\cong\Aut^G(X)$. Denote
by $\pi$ the quotient morphism $X\twoheadrightarrow X/\Gamma$.
\begin{enumerate}
\item  $\X_{G,X/\Gamma}=\{\chi\in \X|
\langle\chi,\Gamma\rangle=1\}$ and $\V_{G,X/\Gamma}=\pi_*(\V)$,
where $\pi_*$ denotes the restriction projection
$\a^*\twoheadrightarrow \a^*_{G,X/\Gamma}$.
\item The map $\pi_*:\D\rightarrow \D_{G,X/\Gamma}, D\mapsto
\pi(D)$,  is the quotient map for the action $\Gamma:\D$. Further,
$G_D=G_{\pi_*(D)},\varphi_{\pi_*(D)}=\pi_*(\varphi_{D})$.
\item If $\Gamma$ is finite, then
$\Lambda_{G,X/\Gamma}=\Lambda$.
\end{enumerate}
\end{Lem}
\begin{proof}
The equality of the weight lattices and assertion 2 are easy. The
equality of the valuation cones follows from Corollary 1.5 in
\cite{Knop5}.  Assertion 3 is a special case of Theorem 6.3 from
\cite{Knop8}. \end{proof}

\begin{Cor}\label{Cor:2.3}
Let $H_1,H_2$ be spherical subgroups of $G$. If
$\X_{G,G/H_1}=\X_{G,G/H_2}$ and $\h_1\sim_G \h_2$, then $H_1\sim_G
H_2$.
\end{Cor}
\begin{proof}
We may assume  that $H_1^\circ=H_2^\circ$. By  Lemma
\ref{Lem:2.1.2}, $\X_{G,G/H_i}=\{\mu\in \X_{G,G/H_i^\circ}| \langle
\mu,H_i/H_i^\circ\rangle=1\}$.  It follows that
$H_1/H_1^\circ=H_2/H_2^\circ$ (as subgroups in
$N_G(H_i^\circ)/H_i^\circ$).
\end{proof}

\subsection{Wonderful embeddings and spherical roots}\label{SUBSECTION_wonderful}
In this subsection $X=G/H$.

 We recall that a $G$-variety
$\overline{X}$ is called a {\it wonderful embedding} of $X$ if
\begin{enumerate}
\item $\overline{X}$ is smooth and projective.
\item There is an open $G$-equivariant embedding $X\hookrightarrow \overline{X}$.
\item $\overline{X}\setminus X$ is a divisor with normal crossings (in particular,
this means that all irreducible components of $\overline{X}\setminus X$ are smooth).
\item Let $D_1,\ldots,D_r$ be irreducible components of $\overline{X}\setminus
X$. Then for any $I\subset \{1,\ldots,r\}$ the subvariety
$\bigcap_{i\in I}D_i\setminus\bigcup_{j\not\in I}D_j$ is a single
$G$-orbit.
\end{enumerate}
It is known that a wonderful embedding of $X$ is unique if it
exists, see, for example, \cite{Luna5} or \cite{Timashev_rev},
Section 30.

\begin{Prop}\label{Prop:2.6}
 $X$ has a wonderful embedding
iff $\X=\Span_\Z(\Psi)$.
\end{Prop}
\begin{proof}
In fact, this is proved in \cite{Knop8}, Corollary 7.2.
\end{proof}

\begin{Cor}\label{Cor:2.1}
Let $\Gamma$ be the annihilator of $\Psi$ in $\A$ and
$\widetilde{H}$ denote  the inverse image of $\Gamma$ in $N_G(H)$.
Then  $G/\widetilde{H}$ admits a wonderful embedding. In particular,
if $H=N_G(H)$, then $X$ admits a wonderful embedding.
\end{Cor}
\begin{proof}
Clearly, $\X_{G,G/\widetilde{H}}=\Span_\Z(\Psi)$. By Lemma
\ref{Lem:2.1.2}, $\Psi_{G,G/\widetilde{H}}=\Psi$. \end{proof}

Thanks to Corollary \ref{Cor:2.1}, all pairs $(\alpha,P),\alpha\in
\Psi,$ correspond to wonderful varieties of rank 1 and so  are
obtained  from those listed in \cite{Wasserman}, Table 1 by a so
called {\it parabolic induction}.

The following proposition describes a relation between $\Psi,\D$. It
was proved by Luna using wonderful varieties.

\begin{Prop}[\cite{Luna4}, Proposition 3.4]\label{Prop:2.4}   For $\alpha\in \Pi(\g)$ exactly one of the following
possibilities holds:
\begin{itemize}
\item[(a)] $\D(\alpha)=\varnothing$.
\item[(b)] $\alpha\in\Psi$. Here $\D(\alpha)=\{D^+,D^-\}$ and
$\varphi_{D^+}+\varphi_{D^-}=\alpha^\vee|_{\a}$.
\item[(c)]  $2\alpha\in \Psi$. In this case $\D(\alpha)=\{D\}$ and $\varphi_{D}=\frac{1}{2}\alpha^\vee|_{\a}$.
\item[(d)] $\Q\alpha\cap \Psi=\varnothing,\D(\alpha)\neq\varnothing$.
Here $\D(\alpha)=\{D\}$  and $\varphi_{D}=\alpha^\vee|_{\a}$.
\end{itemize}
\end{Prop}

We say that $\alpha\in\Pi(\g)$ is of type a) (or b),c),d)) for $X$
if the corresponding possibility of Proposition \ref{Prop:2.4} holds
for $\alpha$.

\begin{Prop}\label{Prop:2}
Let $\alpha\in \Psi,\beta\in \Pi(\g)\cap\Psi, D\in \D(\beta)$. Then
$\langle\varphi_D,\alpha\rangle\leqslant 1$ and the equality holds
iff $\alpha\in \Pi(\g), D\in \D(\alpha)$.
\end{Prop}
\begin{proof}
Using the localization procedure established in \cite{Luna4},
Subsection 3.5, (see also \cite{Luna5}, Subsection 3.2) we reduce
the proof to the case $\Psi=\{\alpha,\beta\}$. In this case
everything follows from the classification in \cite{Wasserman}.
\end{proof}

\subsection{A finiteness result}\label{SUBSECTION_finiteness}
\begin{Prop}\label{Cor:2.2}
Let $Y$ be a locally closed irreducible subvariety of $\Gr_d(\g)$
such that  $\h_0$ is spherical and $\n_\g(\h_0)=\h_0$ for any
$\h_0\in Y$. Then any two elements of $Y$ are $G$-conjugate.
\end{Prop}
\begin{proof}
By \cite{AB}, Corollary 3.2, there is a decomposition
$Y=Y_1\sqcup\ldots\sqcup Y_k$ such that $Y_i$ is an intersection of
$Y$ with a class of $G$-conjugacy of subalgebras. Clearly, $Y_i$ is
 locally closed in $Y$. Let $Y_1$ be open in $Y$ and
$\h_0\in Y_1$. Then  $\h_1\in \overline{G\h_0}$ for any $\h_1\in Y$.
Since $\n_\g(\h_i)=\h_i,i=0,1$, we get $\dim G\h_1=\dim G\h_0$. So
$\h_1\sim_G \h_0$. \end{proof}

\subsection{Inclusions of spherical
subgroups}\label{SUBSECTION_inclusions} Throughout the subsection
$X=G/H$.

The goal of this subsection is to describe the ordered (by
inclusion) set $\H_H$ in terms of $\a,\V,\D$. All results are proved
in \cite{Knop5}, Section 4. We remark that any
$\widetilde{H}\in\H_H$ is spherical. Further,  $\A$ naturally acts
on $\H_H$ preserving the partial order.

\begin{defi}\label{defi:2}
A pair $(\a^1, \D^1)$ consisting of a subspace $\a^1\subset \a^*$
and $\D^1\subset \D$ is called a {\it colored subspace} of
$(\a^*,\D)$ if $\a^1$ coincides with the cone spanned by
$\V\cap\a^1$ and $\varphi_D$ for all $D\in \D^1$. The set of all colored
subspaces of $(\a^*,\D)$ is denoted by $\CS(=\CS_{G,X})$.
\end{defi}

In particular, $(\a^*,\D)\in \CS$. Further, $(\a^1,\varnothing)\in
\CS$  iff $\a^1\subset \V\cap-\V$.

The set $\CS$ has a natural partial order: $(\a^1,\D^1)\preceq
(\a^2,\D^2)$ if $\a^1\subset \a^2, \D^1\subset \D^2$. Further,  $\A$
acts on $\CS$: $a.(\a^1,\D^1)=(\a^1,a.\D^1),a\in \A$.

Choose $\widetilde{H}\in \H_H$. There is the natural morphism
$\pi:X\twoheadrightarrow G/\widetilde{H}, gH\mapsto g\widetilde{H}$.
It induces the inclusion
$\pi^*:\X_{G,G/\widetilde{H}}\hookrightarrow \X$, the projection
$\pi_*:\a^*\twoheadrightarrow \a_{G,G/\widetilde{H}}^*$ and the
embedding $\pi^*:\D_{G,G/\widetilde{H}}\hookrightarrow \D, D\mapsto
\pi^{-1}(D)$. Put $\a_H^{\widetilde{H}}:=\ker\pi_*,
\D_H^{\widetilde{H}}:=\{D\in \D|
\overline{\pi(D)}=G/\widetilde{H}\}=\D\setminus \im\pi^*$. Note that
$\a_H^{\widetilde{H}}$ depends only on
$\a_{G,G/H},\a_{G,G/\widetilde{H}}$.

\begin{Prop}\label{Prop:2.2} $(\a_H^{\widetilde{H}},\D_H^{\widetilde{H}})\in \CS_{G,G/H}$
and  the map $\widetilde{H}\mapsto
(\a_H^{\widetilde{H}},\D_H^{\widetilde{H}})$ is an order-preserving
$\A$-equivariant bijection between $\H_H$ and $\CS$.
\end{Prop}

Now we recover
$\X_{G,G/\widetilde{H}},\V_{G,G/\widetilde{H}},\D_{G,G/\widetilde{H}}$
from $(\a_H^{\widetilde{H}},\D_H^{\widetilde{H}})$.

\begin{Prop}\label{Prop:2.8}
$\X_{G,G/\widetilde{H}}=\{\lambda\in \X| \langle
\lambda,\a_{H}^{\widetilde{H}} \rangle=0\}$,
$\V_{G,G/\widetilde{H}}=\pi_*(\V), G_{\pi^*(D)}=G_D,
\varphi_D=\pi_*(\varphi_{\pi^*(D)})$ for any $D\in
\D_{G,G/\widetilde{H}}$.
\end{Prop}

For example, the subgroup $\widetilde{H}\in \H_H$ is parabolic iff
$\a_{G,G/\widetilde{H}}=\{0\}$ iff $\a_H^{\widetilde{H}}=\a^*$. In
this case $\widetilde{H}\sim_G Q^-$, where $Q^-$ denotes the
parabolic subgroup of $G$ opposite to $Q:=\cap_{D\in
\D\setminus\D_H^{\widetilde{H}}}G_D$ (i.e.,  generated by $B^-$ and
the standard Levi subgroup of $Q$).

\subsection{Around the local structure theorem}\label{SUBSECTION_incl_parab}
Let $Q$ be an arbitrary algebraic group. Choose a Levi decomposition
$Q=M\rightthreetimes \Rad_u(Q)$. Then one can regard $\Rad_u(Q)$ as
a $Q$-variety by setting $m.q=mqm^{-1}, q_1.q=q_1q, m\in M,q,q_1\in
\Rad_u(Q)$. If $X'$ is an $M$-variety, then one regards $X'$ as a
$Q$-variety assuming that $Q$ acts on $X'$ via the projection
$Q\twoheadrightarrow Q/\Rad_u(Q)\cong M$.

\begin{Prop}[The local structure theorem, \cite{Knop3}, Theorem 2.3]\label{Prop:2.5}
Let $X$ be an irreducible smooth  (not necessarily spherical)
$G$-variety and $\D'$ be a finite set of  $B$-stable prime divisors
on $X$. Put $X^0:=X\setminus(\cup_{D\in\D'}D),Q:=\cap_{D\in \D'}G_D$
and let $M$ be the standard Levi subgroup of $Q$. Then there is a
closed $M$-stable subvariety (a {\rm section}) $X'\subset X^0$ such
that the map $\Rad_u(Q)\times X'\rightarrow X^0, (q,x)\mapsto qx,$
is a $Q$-equivariant isomorphism.
\end{Prop}

\begin{Ex}\label{Ex:2.5}
Let $Q^-\in \H_{B^-}$  and $M$ be the standard Levi subgroup of
$Q^-$. Put $Q=BM$. Let $H$ be an algebraic subgroup of $Q^-$ and
$\pi$ denote the natural projection $G/H\twoheadrightarrow G/Q^-$.
Take $\{\pi^{-1}(D), D\in \D_{G,G/Q^-}\}$ for $\D'$. Then one can
take $\pi^{-1}(eQ^-)= Q^-/H$ for $X'$.
\end{Ex}

Below $\D',X^0,Q,M,X'$ have the same meaning as in Proposition
\ref{Prop:2.5}. Note that $X'\cong^M X^0/\Rad_u(Q)$.

\begin{Rem}\label{Rem:2.4.1}
Let $\widetilde{X}$ be another smooth  $G$-variety and
$\pi:\widetilde{X}\rightarrow X$ a smooth surjective $G$-equivariant
morphism with irreducible fibers. Put
$\widetilde{\D}':=\{\pi^{-1}(D), D\in \D'\}$. Then one can take
$\pi^{-1}(X')$ for a section of the action $Q:\pi^{-1}(X^0)$.
\end{Rem}

\begin{Lem}\label{Lem:2.5.1}
If $X$ is affine, then so is $X'$.
\end{Lem}
\begin{proof}
Since $X$ is smooth, any divisor is Cartier. It is known that the complement
to a Cartier divisor in an affine variety is again affine.
So $X^0$ is affine. Being closed in
$X^0$, the variety $X'$ is affine too. \end{proof}

Note that  $X$ is $G$-spherical iff  $X'$ is $M$-spherical.

Until further notices $X$ is a  smooth quasiprojective spherical $G$-variety.

\begin{Lem}\label{Lem:2.3}
\begin{enumerate}
\item
$\X_{M,X'}=\X$. \item The map $\iota:\D_{M,X'}\rightarrow
\D\setminus \D', D\mapsto \Rad_u(Q)\times D$ is bijective and
satisfies $\varphi_{D}=\varphi_{\iota(D)}, M_{D}=M\cap
G_{\iota(D)}$. In particular, $L_{M,X'}=L$.
\item $\Psi_{M,X'}=\Psi\cap \Span_\Q(\Delta(\m))$.
\item $\A\subset \A_{M,X'}$.
\end{enumerate}
\end{Lem}
\begin{proof}
Assertions 1,4 follow from the natural isomorphism between $\K(X')$
and $\K(X)^{\Rad_u(Q)}$. Assertion 2 is straightforward.

We proceed to  assertion 3. The claim  for quasiaffine $X$ follows
from \cite{Losev}, Proposition 8.2. In the general case, since $X$ is quasiprojective,  there is a
quasiaffine $G$-variety $\widetilde{X}$ and a $G$-equivariant
principal $\K^\times$-bundle $\pi:\widetilde{X}\rightarrow X$. Put
$\widetilde{G}:=G\times \K^\times, \widetilde{M}:=M\times
\K^\times$. The group $\widetilde{G}$ acts naturally on
$\widetilde{X}$. Set $\widetilde{\D}'=\{\pi^{-1}(D)| D\in \D'\}$,
$\widetilde{X}':=\pi^{-1}(X')$. Then
$\Psi_{\widetilde{M},\widetilde{X}'}=\Psi_{\widetilde{G},\widetilde{X}}\cap\Span_\Q(\Delta(\m))$.
By \cite{Knop3}, Theorem 5.1,
$\Psi_{M,X'}=\Psi_{\widetilde{M},\widetilde{X}'},
\Psi=\Psi_{\widetilde{G},\widetilde{X}}$. \end{proof}

\begin{Lem}\label{Lem:2.4.2}
\begin{enumerate}
\item Let $G_1$ be a semisimple normal subgroup of $G$. If
$\Pi(\g_1)$ consists of roots of type a), then $G_1$ acts trivially
on $X$.
\item If $\D=\varnothing$, then $(G,G)$ acts trivially on $X$
and $X$ is a single $G/(G,G)$-orbit.\end{enumerate}
\end{Lem}
\begin{proof}
Applying Proposition \ref{Prop:2.5} to $\D'=\D$, we reduce assertion
1 to assertion 2. The latter stems, say, from \cite{Knop3},
Proposition 2.4. \end{proof}

Applying Proposition \ref{Prop:2.5} to $\D'=\D$ and using Lemmas
\ref{Lem:2.3}, \ref{Lem:2.4.2}, we get the following two well-known
lemmas.

\begin{Lem}\label{Prop:2.1}
 $\langle\a_{G,X},\alpha^\vee\rangle=0$ for any simple root
$\alpha$ of type a).
\end{Lem}

\begin{Lem}\label{Prop:2.0}
 $\dim X=\rank \X+\dim G-\dim
P$.
\end{Lem}

Till the end of the subsection $X=G/H$ is a spherical homogeneous
space and $Q^-\in \H_{B^-}\cap\H_H$. Put $X':=Q^-/H\subset G/H$. Let
$M$ denote the standard Levi subgroup of $Q^-$.

\begin{Prop}\label{Prop:2.3}
Let $Q^-,H,M,X'$ be such as above. If $G/H$ admits a wonderful
embedding,  then the following conditions are equivalent:
\begin{itemize}\item[(a)]$\Rad_u(Q^-)\subset H$.\item[(b)]
$\Psi\subset \Span_\Q(\Delta(\m))$.\end{itemize} Under these
equivalent conditions, $\Psi_{M,X'}=\Psi$, $\alpha$ is of type d)
for $X$ and $\D_H^{Q^-}\cap \D(\alpha)=\varnothing$ for any
$\alpha\in \Pi(\g)\setminus \Delta(\m)$, and $
\overline{\Psi}=\overline{\Psi}_{M,X'}$.
\end{Prop}
\begin{proof}
Everything except the equality
$\overline{\Psi}_{M,X'}=\overline{\Psi}$ follows from results of
\cite{Luna5}, Subsection 3.4.  Hence $\D_H^{Q^-}$ is $\A$-stable. By
Proposition \ref{Prop:2.2}, $N_G(H)\subset N_G(Q)=Q$. It follows
that $N_G(H)/H= N_M(H\cap M)/(H\cap M)$ whence
$\Lambda=\Lambda_{M,X'}$ or, equivalently,
$\overline{\Psi}=\overline{\Psi}_{M,X'}$. \end{proof}

There is another distinguished class of inclusions $H\subset Q^-$.
It is well known, see, for example, \cite{Weisfeller}, that  for any
algebraic subgroup  $F\subset G$ there is a parabolic subgroup
$Q\subset G$ such that $F\in \overline{\H}^{Q}$. The following lemma
seems to be standard.
\begin{Lem}\label{Lem:2.4}
Let $X,Q^-,M,X'$ be such as above. Assume, in addition, that $H\in
\overline{\H}^{Q^-}$ and $M\cap H$ is a maximal reductive subgroup
of $H$. Then
 $X'\cong^M M*_{M\cap H}(\Rad_u(\q^-)/\Rad_u(\h))$.
\end{Lem}

\subsection{Some properties of  affine spherical
varieties}\label{SUBSECTION_affine_prelim} 
In this subsection $X$ is a spherical $G$-variety.

\begin{Lem}\label{Lem:2.17}
If $X$ is quasi-affine, then $\varphi_D\neq 0$ for any $D\in \D$.
\end{Lem}
\begin{proof}
There is a regular function vanishing on $D$. By the Lie-Kolchin
theorem there is even a $B$-semiinvariant such function.
\end{proof}

\begin{defi}
By the {\it weight monoid} of $X$ we mean the set
$\X^+_{G,X}:=\{\lambda\in \X(T)| \K[X]^{(B)}_\lambda\neq \{0\}\}$.
\end{defi}

\begin{Lem}\label{Lem:2.16}
\begin{equation}\label{eq:6}\X^+=\{\lambda\in
\X|\langle\lambda,\varphi_D\rangle\geqslant 0, \forall D\in \D\}.
\end{equation}
If $X$ is quasiaffine, then
\begin{equation}\label{eq:5}\X=\Span_\Z(\X^+).\end{equation}
\end{Lem}
\begin{proof}
(\ref{eq:6}) is clear. (\ref{eq:5}) stems from \cite{VP}, Theorem
3.3. \end{proof}

\begin{Prop}\label{Prop:2.18}
Let $X_1,X_2$ be  affine spherical $G$-varieties. Let $X^0_i$ denote
the open $G$-orbit in $X_i,i=1,2$. If $X_1\equiv^G X_2, X^0_1\cong^G
X^0_2$, then $X_1\cong^G X_2$.
\end{Prop}
\begin{proof}
We may consider $\K[X_1],\K[X_2]$ as subalgebras in $\K[X^0_1]\cong
\K[X^0_2]$.  By Lemma \ref{Lem:2.16}, $\X^+_{G,X_1}=\X^+_{G,X_2}$.
The highest weight theory implies  $\K[X_1]=\K[X_2]$.
\end{proof}

\begin{Lem}\label{Lem:2.15}
Let $X_1,X_2$ be  affine spherical $G$-varieties with $X_1\equiv^G
X_2$. Suppose there is a decomposition $G=G^1G^2$ into the locally
direct product such that $X_1\cong^{G^1\times G^2} X_1^1\times
X^2_1$, where $X_1^i$ is a $G^i$-variety. Then $X_1^i$ is a
spherical $G^i$-variety and there is a spherical $G^i$-variety
$X_2^i,i=1,2,$ such that $X_2=X_2^1\times X_2^2$ and
$X_1^i\equiv^{G^i}X_2^i$.
\end{Lem}
\begin{proof}
We may assume that $G=G^1\times G^2$. Everything will  follow if we
check the existence of a $G^i$-variety $X_2^i$ such that
$X_2\cong^{G}X_2^1\times X_2^2$. Put $X^1_2:=X_2\quo G^2,
X_2^2:=X_2\quo G^1$.

By Lemma \ref{Lem:2.16}, $\X_{G,X_1}^+=\X_{G,X_2}^+$. From the
highest weight theory one easily deduces that
$\X_{G^i,X_2^i}^+=\X_{G,X_2}^+\cap \X(T\cap G^i),i=1,2$. Thus
$\X^+_{G^i,X_1^i}=\X^+_{G^i,X_2^i},i=1,2,$ and
$\X^+_{G,X_2}=\X^+_{G^1,X_2^1}+\X^+_{G^2,X_2^2}$. In other words,
$\K[X_2]^{\Rad_u(B)}=\K[X_2^1]^{\Rad_u(B)\cap G^1}\otimes
\K[X_2^2]^{\Rad_u(B)\cap G^2}$. It follows from the highest weight
theory that $\K[X_2]=\K[X_2^1]\otimes \K[X_2^2]$. \end{proof}

\begin{Prop}[\cite{KvS}, Corollary 2.2]\label{Prop:2.10}
If $X$ is smooth and affine, then $X\cong^G G*_HV$, where $H$ is
reductive and $V$ is an $H$-module.
\end{Prop}

 The following lemma seems to be well
known.

\begin{Lem}\label{Lem:2.12}
\begin{enumerate}
\item Let $H_i,i=1,2,$ be a reductive subgroup of $G$ and $V_i$ be
an $H_i$-module. Then $G*_{H_1}V_1\cong^G G*_{H_2}V_2$ iff there are
$g\in G$ and a linear isomorphism $\iota:V_1\rightarrow V_2$ such
that $gH_1g^{-1}=H_2$ and $\iota(hv)=(ghg^{-1})\iota(v)$ for all
$h\in H_1,v\in V_1$.
\item  An element $g\in
N_G(H)/H$ lies in the image of the natural homomorphism
$\Aut^G(X)\rightarrow \Aut^G(G/H)$ iff $V\cong^H V^g$.
\end{enumerate}
\end{Lem}
\begin{proof}
Let $\pi_i:G*_{H_i}V_i\rightarrow G/H_i,i=1,2,$ denote the natural
projection. Let $\varphi:G*_{H_1}V_1\rightarrow G*_{H_2}V_2$ be a
$G$-equivariant isomorphism and $x=[e,0]\in G*_{H_1}V_1$. Clearly,
$G_{\varphi(x)}=H_1$ is conjugate to a subgroup of $H_2$.
Analogously, $H_2$ is conjugate to a subgroup of $H_1$. So we may
assume that $H_1=H_2$. Let $V_2'$ denote the slice module
$T_{\varphi(x)}(G*_{H_1}V_2)/\g_*\varphi(x)$. Note that $\varphi$
gives rise to the linear $H_1$-equivariant isomorphism
$d_x\varphi:V_1\rightarrow V_2'$. Since
$V_2'\cong^{H_1}\pi_2^{-1}(\pi_2(\varphi(x)))$, it follows that
$V_2'\sim^{N_G(H_1)}V_2$ whence the first assertion. Assertion 2
stems easily from the previous argument. \end{proof}

\section{Proofs of Theorems \ref{Thm:1.1},\ref{Thm:1.2}}\label{SECTION_Proof2}
Subsection \ref{SUBSECTION_distinguished} contains  plenty of quite
simple auxiliary results used in the proof of Theorem \ref{Thm:1.2}.
We introduce the notions of a distinguished spherical root
(Definition \ref{defi:1}) and of an automorphism doubling it
(Definition \ref{defi:4}) and show that Theorem \ref{Thm:1.2} is
equivalent to the claim that  any distinguished spherical root has
an automorphism doubling it (Lemma \ref{Lem:4.4}). Other lemmas
study different properties of distinguished roots and automorphisms
doubling them.

To prove our main theorems we use "induction". The "base", roughly, consists of all homogeneous spaces
of the form $G/H$, where $N_G(H)$ does not lie in a proper parabolic subgroup of
$G$. The main goal of  Subsection \ref{SUBSECTION_affine} is to prove "base of induction": assertion 1 of Proposition \ref{Prop:4.1}
for Theorem 2 and Proposition \ref{Prop:4.5} for Theorem 1. The other two results,
assertion 2 of Proposition \ref{Prop:4.1} and Proposition \ref{Prop:4.2}, are of more technical
nature.  Their purposes are described in  the beginning of the subsection.
A common feature of all these results is that they concern smooth affine
spherical varieties.

Subsections \ref{SUBSECTION_Proof2},\ref{SECTION_Proof1} contain "induction
steps" of the proofs of Theorems 2 and 1, respectively. Very roughly,  we need to replace a
homogeneous space $G/H$ with $G/\widetilde{H}$ for an appropriate subgroup $\widetilde{H}\in \underline{\H}_H$.
The proofs are described in more details in  each subsection.

\subsection{Distinguished
roots and automorphisms doubling
them}\label{SUBSECTION_distinguished} Throughout this subsection $X$
is a smooth  spherical $G$-variety.

\begin{defi}\label{defi:1}
An element $\alpha\in \Psi_{G,X}$ is said to  be  {\it
distinguished} if one of the following conditions holds:
\begin{enumerate}
\item $\alpha\in \Pi(\g)$ and  $\varphi_{D}=\frac{1}{2}\alpha^\vee|_{\a_{G,X}}$ for any $D\in\D_{G,X}(\alpha)$.
\item There is a subset $\Sigma\subset \Pi(\g)$ of type
$B_k,k\geqslant 2,$ such that $\alpha=\alpha_1+\ldots+\alpha_k$ and
$\D_{G,X}(\alpha_i)=\varnothing$ for any $i>1$.
\item There is a subset $\Sigma\subset \Pi(\g)$ of type $G_2$
such that $\alpha=2\alpha_2+\alpha_1$.
\end{enumerate}
By the {\it type} of a distinguished root we mean its number in the
previous list. By $\Psi_{G,X}^i, i=1,2,3,$ we denote the subset of
$\Psi_{G,X}$ consisting of all distinguished roots of type $i$. We
put $\widetilde{\alpha}:=\alpha$ (resp.,
$\widetilde{\alpha}=\alpha_1$, $\widetilde{\alpha}=\alpha_2$) for
$\alpha\in \Psi_{G,X}^i,i=1$ (resp., $i=2,3$).  Set
$\Psi^+_{G,X}:=\sqcup_{i=1}^3\Psi_{G,X}^i$.
\end{defi}
Recall that, by our conventions, $\alpha_k,\alpha_2$ denote the
short simple roots in $B_k,G_2$. Note that $\Supp(\alpha)\subset
\{\widetilde{\alpha}\}\cup\Pi(\lfr_{G,X})$ for any $\alpha\in
\Psi^+_{G,X}$.

Distinguished roots were previously considered in Pezzini's paper
\cite{Pezzini1}.

First of all, we study  functorial properties of distinguished roots
(Lemma \ref{Lem:4.1}). For this we need the following auxiliary
lemma that is somewhat similar to Proposition 3.3.2 in \cite{Luna5}.

\begin{Lem}\label{Lem:2.10}
Let $H$ be a spherical subgroup of $G$ and $\widetilde{H}\in \H_H$.
Let $\Psi_0\subset \Psi$ be such that
$\langle\varphi_D,\Psi_0\rangle=0$ for any $D\in
\D_H^{\widetilde{H}}$. Then $\Psi_0\subset \Psi_{G,G/\widetilde{H}}$
and any $\beta\in \Psi_{G,G/\widetilde{H}}\setminus \Psi_0$ is of
the form $\sum_{\alpha\in \Psi\setminus \Psi_0}n_\alpha\alpha,
n_\alpha\geqslant 0$.
\end{Lem}
\begin{proof}
At first, consider the situation when $\a_H^{\widetilde{H}}$ is
generated by $\varphi_D, D\in \D_H^{\widetilde{H}}$. Note that
$\ker\alpha$ is a wall of $\V+\a_H^{\widetilde{H}}$ for any
$\alpha\in \Psi_0$. Thus $\Psi_0\subset \Psi_{G,G/\widetilde{H}}$.
Since $\langle\beta,\V\rangle\leqslant 0$, we obtain
$\beta=\sum_{\alpha\in \Psi}n_\alpha\alpha$ for some $
n_\alpha\geqslant 0$. Put $\beta':=\sum_{\alpha\in
\Psi\setminus\Psi_0}n_\alpha\alpha$. Then
$\langle\beta',\V\rangle\leqslant 0$ and
$\langle\beta-\beta',\a_H^{\widetilde{H}}\rangle=0$. Since
$\ker\beta$ is a wall of $\V+\a_H^{\widetilde{H}}$, we get
$\beta=\beta'$.

We proceed to the general case. We may replace $H$ with the element
of $\H_H$ corresponding to
$(\Span_{\Q}(\varphi_D| D\in \D_{H}^{\widetilde{H}}),\D_{H}^{\widetilde{H}})$. In this case
$\widetilde{H}/H$ is a torus and we use Lemma \ref{Lem:2.1.2}.
\end{proof}

\begin{Lem}\label{Lem:4.1}
\begin{enumerate}
\item Let $X_1,X_2$ be spherical $G$-varieties and $\varphi:X_1\rightarrow
X_2$ be a dominant $G$-equivariant morphism. Then
$\Psi_{G,X_1}^i\cap \Psi_{G,X_2}\subset \Psi_{G,X_2}^i$ for any
$i=1,2,3$. In particular, $\alpha\in \Psi_{G,X_2}^1$ provided
$\alpha\in \Psi^1_{G,X_1}$ and $\#\D_{G,X_2}(\alpha)=2$.
\item Let $H_1$ be a spherical subgroup of $G$, $H_2\in \H_{H_1},H_3\in
\H_{H_2}$. Set $X_l:=G/H_l,l=1,2,3$. If $\alpha\in
\Psi^i_{G,X_1}\cap\Psi^j_{G,X_3}$, then $i=j$ and $\alpha\in
\Psi_{G,X_2}^i$.
\item Let $X_1,X_2$ be spherical $G$-varieties and $\varphi:X_1\rightarrow
X_2$ a generically \'{e}tale $G$-equivariant morphism. Then
$\Psi_{G,X_2}^i= \Psi_{G,X_1}^i\cap\Psi_{G,X_2}$ for any $i$.
\end{enumerate}
\end{Lem}
\begin{proof}
Assertion 1 is clear.
Let us check  assertion 2. 
%
 It follows from  Lemma
\ref{Lem:2.10} that $\alpha\in \Psi_{G,X_2}$. By assertion 1,
$\alpha\in\Psi_{G,X_2}^i, i=j$. Assertion 3 follows easily from the
observation (see \cite{Wasserman}, Table 1) that  for any $\alpha\in
\Psi_{G,X_1}$ either $\alpha$ or $2\alpha$ lies in $\Psi_{G,X_2}$.
\end{proof}

Now let us study the behavior of $\Psi^i_{G,X}$ under some
modifications of the pair $(G,X)$.

\begin{Lem}\label{Lem:4.2}
Let  $\pi:\widetilde{X}\rightarrow X$ be a $G$-equivariant principal
$\K^\times$-bundle, $\widetilde{G}:=G\times\K^\times$. Then
$\Psi^1_{\widetilde{G},\widetilde{X}}\subset \Psi^1_{G,X}$,
$\Psi^i_{\widetilde{G},\widetilde{X}}= \Psi^i_{G,X},i=2,3$.  Let
$\alpha\in \Psi_{G,X}^1$. Then $\alpha\in
\Psi_{\widetilde{G},\widetilde{X}}^1$ iff $\ord_{D_1}(\sigma)=\ord_{D_2}(\sigma)$, where
$\D_{G,X}(\alpha)=\{D_1,D_2\}$ and $\sigma$ is a rational
$B$-semiinvariant section of  $\pi$.
\end{Lem}
\begin{proof}
In the proof we may assume that $X$ is  homogeneous. In this case
$\widetilde{X}$ is a homogeneous $\widetilde{G}$-space. Note that
$\X_{\widetilde{G},\widetilde{X}}=\Z\chi\oplus\X_{G,X}$, where
$\chi$ is the sum of the canonical generator of $\X(\K^\times)$ and
the weight of $\sigma$. Now everything follows from Lemma
\ref{Lem:2.1.2} applied to $\K^\times\subset
\A_{\widetilde{G},\widetilde{X}}$ and Lemma \ref{Lem:4.1}.
\end{proof}

\begin{Lem}\label{Lem:4.3}
Let $\D',Q,M,X'$ be such as in Proposition \ref{Prop:2.5}. Then
$\Psi_{M,X'}^i=\{\alpha\in
\Psi^i_{G,X}|\widetilde{\alpha}\in\Delta(\m)\}, i=1,2,3$.
\end{Lem}
\begin{proof}
By Lemma \ref{Lem:2.3}, $\D_{M,X'}$ is naturally identified with
$\D_{G,X}\setminus \D'$,  $L_{M,X'}=L_{G,X}$, and
$\Psi_{M,X'}=\Psi_{G,X}\cap \Span_\Q(\Delta(\m))$. Since
$\Supp(\alpha)\subset\widetilde{\alpha}\cup\Delta(\lfr_{G,X})$, we
see that the inclusions $\alpha\in \Span_\Q(\Delta(\m))$ and
$\widetilde{\alpha}\in\Delta(\m)$ are equivalent. Now everything
stems from Definition \ref{defi:1}. \end{proof}



\begin{defi}\label{defi:4}
 We say that  $\varphi\in \A_{G,X}$ {\it doubles} $\alpha\in \Psi_{G,X}^+$ if
$\langle\alpha,\varphi\rangle=-1$. The set of all
$\varphi\in\A_{G,X}$ doubling $\alpha$ is denoted by
$\A_{G,X}(\alpha)$.
\end{defi}

Note that $\langle\alpha, \A_{G,X}^\circ\rangle=1$. So
$\A_{G,X}(\alpha)$ is a union of $\A_{G,X}^\circ$-cosets.

\begin{Rem}\label{Rem:4.1} Let $\alpha\in \Psi_{G,X}^1$. The set $\A_{G,X}(\alpha)$ consists
precisely of those $\varphi\in \A_{G,X}$ that transpose elements of
$\D_{G,X}(\alpha)$. This stems from  Proposition \ref{Prop:2.4}
applied to $X$ and $X/\Gamma$, where $\Gamma$ denotes the algebraic
subgroup of $\A_{G,X}$ generated by $\varphi$. If $\alpha\in
\Pi(\g)$ is such that $\D_{G,X}(\alpha)$ contains a $\A$-unstable
divisor, then $\alpha\in \Psi_{G,X}^1$. Conversely, let
$\iota:\D_{G,X}\rightarrow \D_{G,X}$ be a bijection satisfying
$G_D=G_{\iota(D)}, \varphi_D=\varphi_{\iota(D)}$. If $X$ is
homogeneous and Theorem \ref{Thm:1.2} holds for $X$, then $\iota$ is
induced by some element of $\A_{G,X}$.
\end{Rem}

\begin{Lem}\label{Lem:4.4}
 Theorem \ref{Thm:1.2} is equivalent to the following
assertion:
\begin{itemize}
\item[(*)] $\A_{G,X}(\alpha)\neq\varnothing$ for any $\alpha\in \Psi_{G,X}^+$.
\end{itemize}
\end{Lem}
\begin{proof}
We may assume that $X=G/H$. Since
$\overline{\Psi}_{G,X}=\Psi_{G,G/N_G(H)}$,  Theorem \ref{Thm:1.2}
implies (*). From \cite{Wasserman}, Table 1, it follows that
$\overline{\Psi}_{G,X}=\Psi_1\sqcup 2\Psi_2$ for some partition
$\Psi_{G,X}=\Psi_1\sqcup\Psi_2$ with
$\Psi_{G,X}\cap\Lambda(\g)\setminus \Psi^+_{G,X}\subset \Psi_1$. If
(*) holds, then $\Psi_{G,X}^+\subset \Psi_2$. Since the  image  of
$Z(G)$ in $\Aut(X)$ lies in $ \Aut^G(X)$, we see that
$\Psi_{G,X}\setminus \Lambda(\g)\subset \Psi_2$. \end{proof}


\begin{Lem}\label{Lem:4.7}
\begin{enumerate}
\item Let $\alpha\in \Psi_{G,G/H}^+$. Then $\A_{G,G/H^\circ}(\alpha)\neq\varnothing$ is equivalent to $\A_{G,G/H}(\alpha)\neq\varnothing$.
\item  (*) holds for $G/H$ whenever it holds for $G/N_G(H)$.
\end{enumerate}
\end{Lem}
\begin{proof}
Assertion 1 stems from Lemma \ref{Cor:2.4}.

Proceed to assertion 2.  To prove it we need an auxiliary claim.
Namely, let $X=G/H$ be such that $H$ is connected and
$\widetilde{G},\widetilde{X}$ be such as in Lemma \ref{Lem:4.2}. We
want to show that if Theorem \ref{Thm:1.2} holds for $X$, then it
holds for $\widetilde{X}$.

There is a natural exact sequence
$$1\rightarrow \K^\times\rightarrow \A_{\widetilde{G},\widetilde{X}}\rightarrow \A_{G,X}.$$
Let $\mathcal{L}$ be the line $G$-bundle on $X$ corresponding to
$\widetilde{X}\rightarrow X$. The image of the homomorphism
$\A_{\widetilde{G},\widetilde{X}}\rightarrow \A_{G,X}$ coincides
with the stabilizer of $\mathcal{L}$ under the action
$\A_{G,X}:\Pic_G(X)$. Let $\pi$ denote the forgetful  homomorphism
$\Pic_G(X)\rightarrow \Pic(X)$. Since $H$ is connected, we see that
$\X(H)$ is torsion-free. Using the argument of the proof of
\cite{Knop8}, Lemma 7.3, we get
$(\A_{G,X})_{\mathcal{L}}=(\A_{G,X})_{\pi(\mathcal{L})}$. By  Lemma
\ref{Lem:4.2} and  Remark \ref{Rem:4.1}, $\varphi\in
(\A_{G,X})_{\pi(\mathcal{L})}$ iff $\langle\alpha,\varphi\rangle=1$
for all $\alpha\in \Psi^+_{G,X}\setminus
\Psi^+_{\widetilde{G},\widetilde{X}}$. Thus Theorem \ref{Thm:1.2}
holds for $\widetilde{X}$.

Complete the proof of assertion 2. By assertion 1, (*) holds for
$G/N_G(H)^\circ$. Repeatedly applying the auxiliary claim above we
obtain that (*) holds for $X:=G/H^\circ$ considered as a
$\widetilde{G}:=G\times N_G(H)^\circ/H^\circ$-variety. Any $D\in
\D_{G,X}$, any $f\in \K(X)^{(B)}$ and any $v\in \V_{G,X}$ are
$N_G(H)^\circ/H^\circ$-stable (the last two claims follow, for
example, from \cite{Knop4}, Satz 8.1). Moreover,
$\Aut^G(X)=\Aut^{\widetilde{G}}(X)$. From these observations and
Lemma \ref{Lem:2.1} it follows that
$\Psi_{\widetilde{G},X}=\Psi_{G,X},\overline{\Psi}_{\widetilde{G},X}=\overline{\Psi}_{G,X},\Psi^+_{G,X}=
\Psi^+_{\widetilde{G},X}$. Therefore (*) holds for $X$ considered as
a $G$-variety. Applying assertion 1 one more time, we complete the
proof.
\end{proof}

Essentially, the proof of assertion 2 is based on Luna's
augmentation construction, see \cite{Luna5}, Section 6.

\begin{Lem}\label{Lem:4.5}
Let $X=G/H$  and $\alpha\in \Psi_{G,X}^3$. Then
$\A_{G,X}(\alpha)\neq\varnothing$.
\end{Lem}
\begin{proof}
Thanks to Lemma \ref{Lem:4.7}, we may assume that $G/H$ admits a
wonderful embedding. Note that $G=G'\times G_2$ for some reductive
group $G'$. In the notation of Definition \ref{defi:1},
$\langle\alpha_1^\vee,\Psi_{G,X}\rangle=0$. Hence
$\Psi_{G,X}=\{\alpha\}\sqcup (\Psi_{G,X}\cap\Span_\Q(\Delta(\g')))$.
It follows from  Propositions \ref{Prop:2.4},\ref{Prop:2} that
$\langle\varphi_D,\alpha\rangle=0$ for any $D\in \D_{G,X}\setminus
\D_{G,X}(\widetilde{\alpha})$. It follows from \cite{Luna5},
Proposition 3.5, that $H=H_0\times H'$, where $H_0\subset G_2,
H'\subset G'$. So we may assume that $G=G_2$. In this case $H=A_2$
and the claim follows.
\end{proof}

\subsection{Case of smooth affine $X$}\label{SUBSECTION_affine}
This subsection studies some properties of smooth
affine spherical varieties. Proposition \ref{Prop:4.1} deals with
some special (in fact, almost all) spherical affine homogeneous
spaces.  Assertion 1  proves Theorem \ref{Thm:1.2} for these
homogeneous spaces, while assertion 2 is an auxiliary result used in
the proof of Theorem \ref{Thm:1.1}.
Two assertions of Proposition
\ref{Prop:4.1} are brought together because their proofs use similar ideas.
 Proposition \ref{Prop:4.5}
proves Theorem \ref{Thm:1.1} for some special reductive subgroups of
$G$. Proposition \ref{Prop:4.2} establishes a relation between the
sets of distinguished roots for an affine homogeneous vector bundle
and its closed orbit. This result is very important in the proof of
Proposition \ref{Prop:4.4} in Subsection \ref{SUBSECTION_Proof2}.
\begin{Prop}\label{Prop:4.1}
Let $H$ be a reductive subgroup of $G$ such that $N_G(H)$ is not
contained in a proper parabolic subgroup of $G$. Then
\begin{enumerate}
\item The condition (*) of Lemma \ref{Lem:4.4} holds for $X=G/H$.
\item $g^2\in N_G(H)^\circ Z(G)$ for any $g\in N_G(H)$.
\end{enumerate}
\end{Prop}
\begin{Lem}\label{Lem:4.6}
Let $X$ be an arbitrary smooth spherical $G$-variety. Then there is
an epimorphism $\Pic(X)\twoheadrightarrow
\Z^{\Psi_{G,X}^1}$.
\end{Lem}
\begin{proof}
This stems easily from the observation that $\Pic(X)\cong
\Z^{\D_{G,X}}/((f_\lambda), \lambda\in \X_{G,X})$, see \cite{Brion4}
or Section 17 of \cite{Timashev_rev}.
\end{proof}

\begin{proof}[Proof of Proposition \ref{Prop:4.1}]
It is enough to check both assertions only when $H=N_G(H)^\circ$.
This is clear for assertion 2 and follows from Lemma \ref{Lem:4.7}
for assertion 1. One may assume, in addition, that $\h$ is
indecomposable, i.e., there are no proper complementary ideals
$\g_1,\g_2\subset \g$ such that $\h=(\h\cap\g_1)\oplus(\h\cap\g_2)$.
In particular, $G$ is semisimple.
 If $\rank_G(X)=1$,
then the pair $(\g,\h)$ is presented in Table 1 of \cite{Wasserman}
and we use a case by case argument. Below we assume that
$\rank_G(X)>1$.

Note that if $H\neq (H,H)$, then  $(H,H)$ is not spherical and
$N_G(H)$ acts  on $\z(\h)$ without nonzero fixed vectors. The latter
follows from the assumption that $N_G(H)$ is not contained in a
parabolic. If $(H,H)$ is spherical, then the group
$N_G((H,H))/(H,H)$ is abelian whence $N_G(H)$ and $Z(H)^\circ$
commute.

Connected  reductive spherical subgroups $H\subset G$ with
$N_G(H)^\circ=H$  were classified in \cite{Kramer} ($G$ is simple)
and \cite{Brion_class},\cite{Mikityuk} ($G$ is not simple).  In the
case of simple  $G$  the  monoid $\X^+_{G,G/H}$ is presented in
Tabelle of \cite{Kramer}. Recall that
$\X_{G,X}=\Span_\Z(\X^+_{G,X})$, see Lemma \ref{Lem:2.16}.

Let us prove assertion 1. We may assume that $G$ is of adjoint type.
 Since $N_G(H)$ acts  on $\z(\h)$ without nonzero fixed vectors, we
get $\#\Pic(G/N_G(H))<\infty$. By Lemma \ref{Lem:4.6},
$\Psi^1_{G,G/N_G(H)}=\varnothing$.  In other words,
$\A_{G,X}(\alpha)\neq\varnothing$ for any $\alpha\in \Psi_{G,X}^1$.

Now let $\alpha\in \Psi_{G,X}^2$ and $\g_1$ be a unique simple ideal
of $\g$ with $\alpha\in \Pi(\g_1)$. Then $\g_1\cong
B_n,C_n,n>1,F_4,$ and $\langle\a_{G,X},\beta^\vee\rangle=0$ for some
short simple root $\beta\in \Pi(\g_1)$. It follows from
\cite{Kramer},\cite{Brion_class},\cite{Mikityuk} that $\g$ does not
have an ideal of type $F_4$. Suppose $\g_1\cong B_n,n>2$. Using the
computation of $\a_{G,X}$ (\cite{ranks}), we see that $G=\SO_{2n+1},
H=\SO_{n_1}\times\SO_{2n+1-n_1},1<n_1< [\frac{n}{2}]$.  But here
$\D_{G,X}(\alpha_i)=\varnothing, i>n_1,$ and
$\alpha_{n_1}+\ldots+\alpha_n\not\in \X_{G,X}$. So
$\g_1\cong\sp_{2k}$ and $\D_{G,X}(\alpha_{k-1})=\varnothing,
\alpha_k+\alpha_{k-1}\in \Psi_{G,X}$, where $\alpha_{k-1},\alpha_k$
are simple roots in $C_k$. Since $H$ is semisimple or $(H,H)$ is not
spherical, we get from \cite{Brion_class}, \cite{Mikityuk} that all
summands of $\g$ are of type $C$. Using computations of $\a_{G,X}$
from \cite{ranks}, we obtain
$\g=\sp_{4n},\h=\sp_{2n}\times\sp_{2n}$. In this case indeed
$\#\Psi_{G,X}^2=1$  and a nontrivial element of $N_G(H)/H$ doubles
the element from $\Psi^2_{G,X}$.

Proceed to assertion 2. At first, suppose  $G/H$ is a symmetric
space. We may assume that $G$ is simply connected. Then $H=G^\sigma$
for some involutory automorphism $\sigma$ of $G$, see \cite{VO}. It
follows from \cite{Vust}, Lemma 1 in 2.2, that for any $g\in N_G(H)$
there is $z\in Z(G)$ such that $\sigma(g)=gz$. Assertion 2 will
follow if we check $\sigma(g^2z^{-1})=g^2z^{-1}$. Since
$\sigma^2=id$, we get $g=\sigma(gz)=g\sigma(z)z$ whence
$\sigma(z)=z^{-1}$. So
$\sigma(g^2z^{-1})=g^2\sigma(z)^2\sigma(z^{-1})=g^2\sigma(z)=g^2z^{-1}$.

Consider the situation when $G/H$ is not symmetric. Here we assume
that $G$ is of adjoint type. One checks that if $G$ is not simple,
then $N_G(H)=H$. So it remains to consider the following pairs
$(\g,\h)$: $(\so_{2n+1},\gl_n), (\sp_{2n},\sp_{2n-2}\times\so_2),
(\so_{10},\spin_7\times\so_2),(\so_9,\spin_7),(\so_8,G_2)$. One
checks case by case that $N_G(H)/H\cong \Z_2$ or $1$.
\end{proof}

\begin{Prop}\label{Prop:4.5}
Let $H_1,H_2$ be reductive algebraic subgroups of $G$. Suppose
$H_i,i=1,2,$ is not  contained in any proper parabolic subgroup of
$G$ and $H_i=N_G(H_i)$. If $H_1\equiv^G H_2$, then $H_1\sim_G H_2$.
\end{Prop}
\begin{proof}
 At first, consider the case of simple $G$. Since
$\a_{G,G/H_1}=\a_{G,G/H_2}$, inspecting Tabelle in \cite{Kramer}, we
see that $\g=\so_{2n+1}, \h_1=\gl_n, \h_2=\so_n\times\so_{n+1}$. But
$$\X_{\SO_{2n+1},\SO_{2n+1}/\GL_n}=\Span_\Z(\pi_1,\ldots,\pi_{n-1},2\pi_n),$$
$$\X_{\SO_{2n+1},\SO_{2n+1}/\SO_n\times\SO_{n+1}}=\Span_\Z(2\pi_1,\ldots,2\pi_n),$$
and $\X_{G,G/H_1}=\X_{G,G/H_2}$ is a subgroup of index 2 in both
these lattices, which is absurd.

Let us consider the case when $G$ is not simple. Thanks to Lemma
\ref{Lem:2.15}, we may assume that both $\h_1,\h_2$ are
indecomposable.

Suppose $G/H_1$ is a symmetric space. Then we may assume that
$G\cong H_1\times H_1$ and $H_1$ is embedded diagonally into $G$. It
is clear that $\X_{G,G/H_i}\cap \Span_\Q(\Pi(\g_1))=\{0\}$ for any
simple ideal $\g_1\subset\g$. It follows that the projection of
$\h_2$ to both ideals of $\g$ is surjective. Thence $H_1\cong H_2$.
From $\X_{G,G/H_1}=\X_{G,G/H_2}$ we deduce that $H_1\sim_G H_2$.

It remains to consider the case when neither $G/H_1$ nor $G/H_2$ is
symmetric. Using the classification of
\cite{Brion_class},\cite{Mikityuk} and the equality $\dim H_1=\dim
H_2$, we see that $G=\Sp_{4}\times\Sp_{2m}\times\Sp_{2n},$ $
H_1,H_2\cong\SL_2\times\SL_2\times \Sp_{2m-2}\times\Sp_{2n-2}$,
where one of the normal subgroups $\SL_2\subset H_1$ is embedded
into $\Sp_4$, while the normal subgroups $\SL_2\subset H_2$ are
embedded diagonally into $\Sp_4\times\Sp_{2n},\Sp_4\times\Sp_{2m}$.
By the Frobenius reciprocity, $\pi_1'+\pi_1''\in
\X^+_{G,G/H_1}\setminus \X^+_{G,G/H_2}$ (here $\pi_1',\pi_1''$
denote the highest weights of tautological representations of
$\Sp_{2m},\Sp_{2n}$, respectively). Contradiction with (\ref{eq:6}).
 \end{proof}

\begin{Prop}\label{Prop:4.2}
Let $X$ be a spherical $G$-variety of the form $G*_HV$, where $H$ is
a reductive subgroup of $G$ and $V$ is an $H$-module. Then
$\Psi_{G,X}^i\subset\Psi_{G,G/H}^i$ for  $i=1,2$.
\end{Prop}
\begin{proof}
Let $\pi$ denote the natural projection $G*_HV\twoheadrightarrow
G/H$ and $\pi^*:\D_{G,G/H}\hookrightarrow \D_{G,G*_HV}$ be the
corresponding embedding.

{\it Step 1.} Let us show that if $H=G$, then
$\Psi^i_{G,X}=\varnothing$. For $i=1$ this follows from Lemma
\ref{Lem:4.6}, since $X$ is factorial. Now suppose  $\alpha\in
\Psi_{G,X}^2$. In the notation of Definition \ref{defi:1},
$\D_{G,X}(\alpha_i)=\varnothing, i=\overline{2,k},
\D_{G,X}(\alpha_1)\neq\varnothing$. In other words, any  element in
$\K[X]^{(B)}$ is $P_{\alpha_{i}}$-semiinvariant for $i=2,\ldots,k$,
but there is a prime element in $\K[X]^{(B)}$, that is not
$P_{\alpha_1}$-semiinvariant. From the classification of spherical
$G$-modules (see, for example, \cite{Leahy}, Theorems 1.4,2.5 and
Tables 1,2, or \cite{Knop11}, Section 5) one gets
$G=\SO(2k+1)\otimes \K^\times$ as a linear group. However, here
$\Psi_{G,V}=\{2(\alpha_1+\ldots+\alpha_k)\}$.

{\it Step 2.} Let us show that $\Psi_{G,X}^i=\varnothing, i=1,2,$
whenever $(G,G)\subset H$. Indeed, replacing $G$ with a covering, we
may assume that $G=Z\times H^\circ$, where $Z\cong (\K^\times)^m$.
Now there is an \'{e}tale $G$-equivariant morphism $Z\times
V\rightarrow X$, where $G$ acts on $V$ via the projection
$G\twoheadrightarrow H^\circ$. By assertion 3 of Lemma
\ref{Lem:4.1}, $\Psi_{G,X}^i\subset\Psi_{G,Z\times V}^i$. Obviously,
$\Psi^i_{G,Z\times V}=\Psi^i_{H^\circ,V}$.  By the previous step,
$\Psi_{G,Z\times V}^i=\varnothing$.

{\it Step 3.} Let us check that
$\D_{G,X}(\widetilde{\alpha})=\pi^*(\D_{G,G/H}(\widetilde{\alpha}))$.
Assume the contrary: $\overline{\pi(D)}=G/H$  for some $D\in
\D_{G,X}(\widetilde{\alpha})$. Using  Lemma \ref{Lem:2.17}, we see
that
 $\overline{\pi(D)}=G/H$ for {\it any} $D\in
 \D_{G,X}(\widetilde{\alpha})$. Thus
$\D_{G,G/H}(\widetilde{\alpha})=\varnothing$.  Applying Proposition
\ref{Prop:2.5} to  $(X,\pi^*(\D_{G,G/H}))$, $(G/H,\D_{G,G/H})$, we
get the Levi subgroup $M=L_{G,G/H}$ and sections $X_0'\subset
G/H,X':=\pi^{-1}(X_0')$. By Lemma \ref{Lem:4.3}, $\alpha\in
\Psi_{M,X'}^i$. Note that $\D_{M,X'_0}=\varnothing$. From Lemma
\ref{Lem:2.4.2} it follows that $X'\cong^M M*_S V$, where
$(M,M)\subset S$.  By step 2, $\Psi_{M,X'}^i=\varnothing$.
Contradiction.

{\it Step 4.} If $\alpha\in \Psi_{G,X}^1$, then everything follows
from step 3. So assume that $\alpha\in \Psi_{G,X}^2$. Let
$X^0=G/H^0$ denote the open $G$-orbit in $X$. Since fibers of the
natural morphism $G/H^0\rightarrow G/H$ are irreducible, we see that
$H\in \H_{H^0}$. Since
$\#\D_{G,G/H}(\widetilde{\alpha})=\#\D_{G,G/H^0}(\widetilde{\alpha})=1$,
it follows from Propositions \ref{Prop:2.4},\ref{Prop:2} that
$\langle\varphi_D,\alpha\rangle\leqslant 0$ for any $D\in
\D_H^{H^0}$. Therefore $\langle\varphi_D,\alpha\rangle=0$ for all
such $D$. It remains to use Lemma \ref{Lem:2.10} and the first
assertion of Lemma \ref{Lem:4.1}.
\end{proof}
\subsection{Completing the proof of Theorem \ref{Thm:1.2}}\label{SUBSECTION_Proof2}
The scheme of the proof is as follows.  First, using results of Subsections \ref{SUBSECTION_distinguished},\ref{SUBSECTION_affine},
we explain how to reduce the proof to the case when $N_G(H)=H$
and $H$ is contained in a proper parabolic. Here we need
to show that there are  no distinguished spherical roots. Propositions \ref{Prop:4.4},
\ref{Prop:4.3} together contradict the existence of distinguished roots. More precisely,
if there is a distinguished root
$\alpha\in \Psi_{G,G/H}$ then we show that
there exists $\widetilde{H}\in \underline{\H}_H$ with $\alpha\in \Psi_{G,G/\widetilde{H}}$
(Proposition \ref{Prop:4.4}). But, according to Proposition \ref{Prop:4.3}, existence of such $\widetilde{H}$ contradicts the inductive assumption  in the next paragraph.

So let $X=G/H$ be a spherical homogeneous space. In the proof we may
assume that Theorem \ref{Thm:1.2} is proved for all groups $G_0$
with $\dim G_0<\dim G$ and for all spherical homogeneous spaces
$G/H_0$ with $\dim G/H_0<\dim G/H$. By  Lemma \ref{Lem:4.7} and
Proposition \ref{Prop:2.3}, we may assume that $H=N_G(H)$ and
$\Rad_u(Q)\not\subset H$ whenever $H\subset Q$ for some parabolic
subgroup $Q\subset G$.

By assertion (1) of Proposition \ref{Prop:4.1}, $H$ is contained in some proper parabolic subgroup of
$G$. Thanks to Lemma \ref{Lem:4.4}, Theorem \ref{Thm:1.2} stems from the following two
propositions.

\begin{Prop}\label{Prop:4.4}
Let $\alpha\in \Psi_{G,X}^i$  for some $i=1,2$. Suppose $H$ is
contained in some proper parabolic subgroup of $G$. Then there is
$\widetilde{H}\in \underline{\H}_H$ such that $\alpha\in
\Psi_{G,G/\widetilde{H}}$.
\end{Prop}

\begin{Prop}\label{Prop:4.3}
Let $\alpha\in \Psi_{G,X}^i,i=1,2$. There is no $\widetilde{H}\in
\underline{\H}_H$ such that $\alpha\in \Psi_{G,G/\widetilde{H}}$.
\end{Prop}

\begin{Lem}\label{Lem:4.8}
Let $\widetilde{H}$ be a spherical subgroup of $G$, $\alpha\in
\Psi_{G,G/\widetilde{H}}^i, i=1,2$, $Q^-\in \H_{B^-}$, and $M$ be
the standard Levi subgroup of $Q^-$.
\begin{enumerate}
\item If $\widetilde{H}\subset Q^-$ and $i=2$, then $\widetilde{\alpha}\in \Delta(\m)$.
\item If $i=1, N_G(\widetilde{H})\subset Q^-$, and $\A_{G,G/\widetilde{H}}(\alpha)\neq\varnothing$,
then  $\alpha\in \Delta(\m)$.
\end{enumerate}
\end{Lem}
\begin{proof}
The inclusions $\widetilde{\alpha}\in \Delta(\m)$ and
$\D_{G,G/\widetilde{H}}(\widetilde{\alpha})\subset \D_{\widetilde{H}}^{Q^-}$ are
equivalent.  Note that, thanks to Propositions
\ref{Prop:2.4},\ref{Prop:2},
$\langle\varphi_D,\alpha\rangle\leqslant 0$ whenever $D\not\in
\D_{G,G/\widetilde{H}}(\widetilde{\alpha})$. Therefore
$\D_{\widetilde{H}}^{Q^-}\cap
\D_{G,G/\widetilde{H}}(\widetilde{\alpha})\neq\varnothing$. This observation
proves assertion 1. To prove assertion 2 we note that
$\D_{\widetilde{H}}^{Q^-}$ is $\A_{G,G/\widetilde{H}}$-stable whence
$\D_{\widetilde{H}}^{Q^-}\supset
\D_{G,G/\widetilde{H}}(\widetilde{\alpha})$.
\end{proof}

\begin{Lem}\label{Lem:5}
\begin{enumerate}
\item
If $\widetilde{H}\in \underline{\H}_H$, then
$[\Rad_u(\widetilde{\h}),\Rad_u(\widetilde{\h})]\subset \Rad_u(\h)$.
\item If $H$ is contained in a parabolic subgroup $Q\subset G$ and
$\Rad_u(Q)\not\subset H$, then
$\underline{\H}_H\cap\H^{H\Rad_u(Q)}\neq \varnothing$.
\end{enumerate}
\end{Lem}
\begin{proof}
Note that $\Rad_u(\q)$ is a nilpotent Lie algebra.

Put $\n:=\Rad_u(\h), \widetilde{\n}:=\Rad_u(\widetilde{\h})$. The
projection of $\n$ to
$\widetilde{\n}/[\widetilde{\n},\widetilde{\n}]$ is not surjective.
Otherwise $\n$ generates $\widetilde{\n}$ whence coincides with
$\widetilde{\n}$. Since $\widetilde{\n}/\n$ is an irreducible
$H$-module, we see that $[\widetilde{\n},\widetilde{\n}]\subset \n$.

Proceed to assertion 2. Let $\q_i, i=0,1,2,\ldots$ be the sequence
of ideals of $\q$ defined inductively by $\q_0=\Rad_u(\q),
\q_i=[\q_0,\q_{i-1}]$. Choose the minimal number $i$ such that
$\q_i\subset \h$. Let $\v$ be an irreducible $H$-submodule in
$\q_{i-1}/(\q_{i-1}\cap\h)$, $\widetilde{\n}$ be the inverse image
of $\v$ in $\q_{i-1}$, and $\widetilde{N}$ the connected subgroup of
$\Rad_u(Q)$ corresponding to $\widetilde{\n}$. Then
$H\widetilde{N}\in \underline{\H}_H\cap \H^{H\Rad_u(Q)}$.
\end{proof}

\begin{proof}[Proof of Proposition \ref{Prop:4.4}]
Assume the contrary.

{\it Step 1.} Let us show that there exists $Q^-\in \H_{B^-}$ such
that:
\begin{enumerate}
\item $\widetilde{\alpha}\in \Delta(\m)$, where $M$ is the standard Levi subgroup of
$Q^-$.
\item There is $g\in G$ such that $gHg^{-1}\in \overline{\H}^{Q^-}$.
\end{enumerate}

If $\alpha\in \Psi_{G,X}^2$, then the claim follows from Lemma
\ref{Lem:4.8}. So we may assume that $\alpha\in \Psi_{G,X}^1$. Let
us check that there is a $G$-equivariant principal
$\K^\times$-bundle $X_1\rightarrow X$ such that $X_1$ is quasiaffine
and $\alpha\in \Psi_{G_1,X_1}^1,G_1:=G\times \K^\times$.  For $X_1$
we take the principal $\K^\times$-bundle over $X$ corresponding to a
 divisor of the form
$m\sum_{D\in \D_{G,X}} D$. For sufficiently large $m$ the divisor
$m\sum_{D\in \D_{G,X}} D$ is very ample on the wonderful embedding
$\overline{X}$ of $X$ (see \cite{Brion4}, Section 2). Being an open
subset in the affine cone over $\overline{X}$, the variety $X_1$ is
quasiaffine. By Lemma \ref{Lem:4.2}, $\alpha\in \Psi^1_{G_1,X_1}$.
Clearly, $X_1$ is a homogeneous $G_1$-space. Let $H_1$ denote the
stabilizer of a point in $X_1$ lying over $eH$. So the projection of
$H_1$ to $G$ coincides with $H$. By Sukhanov's theorem,
\cite{Sukhanov}, there are a parabolic subgroup $Q^-\subset G$,  a
$G_1$-module $V$, and a nonzero $Q_1^-$-semiinvariant vector $v\in
V$ such that $H_1\in \overline{\H}^{(Q^-_1)_v}$, where
$Q^-_1:=Q^-\times\K^\times$. We may assume that $B^-\subset Q^-$. By
assertion 2 of Lemma \ref{Lem:5}, there is $\widetilde{H}_1\in
\underline{\H}_{H_1}\cap \overline{\H}^{Q^-_1}$. Automatically,
$\widetilde{H}_1\in \overline{\H}^{(Q_1^-)_v}$. From Sukhanov's
theorem it follows that $G_1/\widetilde{H}_1$ is quasiaffine.   By
Lemma \ref{Lem:2.17}, $\varphi_D\neq 0$ for all $D\in \D_{G_1,G_1/\widetilde{H}_1}$. Therefore
$\#\D_{G_1,X_1}(\alpha)\cap\D_{H_1}^{\widetilde{H}_1}\neq 1$, in
other words,  $\#\D_{G_1,G_1/\widetilde{H}_1}(\alpha)=0$ or $2$.

Let $\widetilde{H}$ denote the projection of $\widetilde{H}_1$ to
$G$. Obviously, $G_1/\widetilde{H}_1\rightarrow G/\widetilde{H}$ is
a $G$-equivariant principal $\K^\times$-bundle. It follows that
$\#\D_{G_1,G_1/\widetilde{H}_1}(\alpha)=
\#\D_{G,G/\widetilde{H}}(\alpha)$. By our assumptions,
$\D_{G,G/\widetilde{H}}(\alpha)=\varnothing$. Since
$\widetilde{H}\subset Q^-$, we get
$\D_{G,G/Q^-}(\alpha)=\varnothing$, q.e.d.

{\it Step 2.} Let $Q^-,M$ satisfy (1),(2) with $g=1$. We may assume
that  $S:=H\cap M$ is a maximal reductive subgroup of $H$. Note
that, by our assumptions, $H\subsetneq \widehat{H}:=H\Rad_u(Q)$. Put
$\underline{X}:=G/\widehat{H}, X'=Q^-/H\cong^M
M*_S(\Rad_u(\q^-)/\Rad_u(\h)),
\underline{X}':=Q^-/\widehat{H}\cong^M M/S$. By Lemma \ref{Lem:4.3},
$\alpha\in \Psi_{M,X'}^i$. Thanks to Proposition \ref{Prop:4.2},
$\alpha\in \Psi_{M,\underline{X}'}^i$. Applying Lemma \ref{Lem:4.3}
again, we see that $\alpha\in \Psi_{G,\underline{X}}^i$. By Lemma
\ref{Lem:5}, there is $\widetilde{H}\in\underline{\H}_H\cap
\H^{\widehat{H}}$. Assertion 2 of Lemma \ref{Lem:4.1} implies
$\alpha\in \Psi_{G,G/\widetilde{H}}^i$. \end{proof}

\begin{proof}[Proof of Proposition \ref{Prop:4.3}]
Assume the contrary: let such $\widetilde{H}$ exist.  We may assume
that there is $Q^-\in \H_{B^-}$ such that $N_G(\widetilde{H})\in
\overline{\H}^{Q^-}$ and $M\cap N_G(\widetilde{H})$ is a maximal
reductive subgroup of $N_G(\widetilde{H})$, where, as usual, $M$ is
the standard Levi subgroup of $Q^-$. Put $S:=M\cap H$. It is a
maximal reductive subgroup of $\widetilde{H}$ or, equivalently, of
$H$. Set $\underline{X}:=G/\widetilde{H}$ and let $\rho$ denote the
natural projection $X\twoheadrightarrow \underline{X}$. Put
$X':=Q^-/H, \underline{X}':=Q^-/\widetilde{H}$.

Let $\pi_\alpha$ be an element of $\a_{G,X}^*$ such that
$\langle\pi_\alpha,\beta\rangle=\delta_{\alpha\beta}$ for all
$\beta\in \Psi_{G,X}$. Define $\varphi_\alpha\in A_{G,X}$ by
$\langle\varphi_{\alpha},\lambda\rangle=\exp(\pi i
\langle\pi_\alpha, \lambda\rangle)$ (actually, $\varphi_\alpha\in
\Hom_\Z(\X_{G,X},\overline{\Q}^\times)$, where $\overline{\Q}$
denotes the algebraic closure of $\Q$).

{\it Step 1.} By our assumptions, Theorem \ref{Thm:1.2} holds for
$X',\underline{X}$. Applying Lemma \ref{Lem:4.8} to $\underline{X}$,
we get $\widetilde{\alpha}\in \Delta(\m)$. By Lemma \ref{Lem:4.3},
$\alpha\in \Psi_{M,X'}^i$.

Let us check that $\varphi_\alpha\in \A_{M,X'}$. Since Theorem
\ref{Thm:1.2} holds for $X'$, we get
$\overline{\Psi}_{M,X'}=\{2\alpha\}\sqcup \Psi_1$, where
$\Psi_1\subset\Span_\Q(\Psi_{G,X}\setminus\{\alpha\})$. Thus
$\langle\varphi_\alpha,\overline{\Psi}_{M,X'}\rangle=1$, q.e.d.

Let us check that $\rho_*(\varphi_\alpha)\in \A_{G,\underline{X}}$,
where $\rho_*$ denotes the homomorphism $A_{G,X}\rightarrow
A_{G,\underline{X}}$ induced by $\rho$.  By Lemma \ref{Lem:2.10},
$\overline{\Psi}_{G,\underline{X}}=\{2\alpha\}\sqcup \Psi_2$, where
$\Psi_2\subset \Span_\Q(\Psi_{G,X}\setminus \{\alpha\})$. Again,
$\langle\rho_*(\varphi_\alpha),\overline{\Psi}_{G,\underline{X}}\rangle=1$.

{\it Step 2.} Let $\gamma$ denote the image of $\varphi_\alpha$ in
$\Aut^M(M/S)$ (see Lemma \ref{Lem:2.11}). Put
$\widetilde{\n}:=\Rad_u(\widetilde{\h}), \n:=\Rad_u(\h)$,
$\v:=\widetilde{\h}/\h$. Clearly, $\v$ is an irreducible $S$-module.
Let us check that $\v\cong^S \v^{\gamma}$. Thanks to step 1,
$\varphi_\alpha\in \Aut^M(X')$. By Lemma \ref{Lem:2.2} and assertion
4 of Lemma \ref{Lem:2.3}, $\rho_*(\varphi_\alpha)\in
\Aut^M(\underline{X}')$. Lemma \ref{Lem:2.12} implies that
$\Rad_u(\q^-)/\n\cong^S (\Rad_u(\q^-)/\n)^\gamma,
\Rad_u(\q^-)/\widetilde{\n}\cong^S
(\Rad_u(\q^-)/\widetilde{\n})^\gamma$ whence $\v\cong^S\v^\gamma$.

{\it Step 3.} Put $\Gamma:=(M\cap N_G(\widetilde{H}))/S$. From the
assumptions on $Q^-,M$ it follows that the natural homomorphism
$\Gamma\rightarrow N_G(\widetilde{H})/\widetilde{H}$ is an
isomorphism. The inverse of this isomorphism is induced by the
projection $Q^-\twoheadrightarrow Q^-/\Rad_u(Q^-)\cong M$. Under the
identification $\Gamma\cong \A_{G,\underline{X}}$, we have
$\gamma=\rho_*(\varphi_\alpha)$. Let $Y$ denote the subvariety in
$\Gr_{\dim\n}(\widetilde{\n})$ consisting of all $S$-stable
subspaces  $\n_0$ such that $[\widetilde{\n},\widetilde{\n}]\subset
\n_0$ and $\widetilde{\n}/\n_0$ is an irreducible $S$-module.
Clearly, $Y$ is a disjoint union of projective spaces. Let $Y_\v$
denote the component of $Y$ consisting of all $\n_0\in Y$ such that
$\widetilde{\n}/\n_0\cong^S \v$. By assertion 1 of Lemma
\ref{Lem:5}, $\n\in Y_\v$. There is a natural action $\Gamma:Y$. The
subgroup $\Gamma_\v\subset \Gamma$ (the stabilizer of $\v$)
coincides with the stabilizer of $Y_\v$ under this action.    The
subset $Y^0_\v:=\{\n_0\in Y_\v| \n_\g(\s+\n_0)=\s+\n_0\}$ is open
and $\Gamma_\v$-stable. By Proposition \ref{Cor:2.2}, $\s+\n_0\sim_G
\h$ for any $\n_0\in Y_\v^0$.

Let $\n_1,\n_2\in Y$ be such that $g(\s+\n_1)=(\s+\n_2)$ for some
$g\in \widetilde{H}$. It is easily seen that
$\n_1=\widetilde{\n}\cap g(\s+\n_1)$ whence $\n_1=\n_2$. So
$\n_1,\n_2\in Y$ are $\Gamma$-conjugate iff $\s+\n_1,\s+\n_2$ are
$N_G(\widetilde{H})$-conjugate.

Let us check that $N_G(\widetilde{H})$ acts on
$\H^{\widetilde{H}}\cap\{gHg^{-1}, g\in G\}$ with finitely many
orbits. Set $A:=\{g\in G| gHg^{-1}\in \H^{\widetilde{H}}\}=\{g\in G|
g^{-1}\widetilde{H}g\in \H_{H}\}$.  Clearly, $A$ is a union of left $N_G(\widetilde{H})$-cosets. 
One easily deduces from Proposition \ref{Prop:2.2} that the subset
$\{g^{-1}\widetilde{H}g,g\in A\}\subset \H_H$ is finite. So $A$
consists of finitely many $N_G(\widetilde{H})$-cosets.

So $Y^0_\v$ consists of finitely many $\Gamma_\v$-orbits. There is a
natural inclusion $\Gamma_\n\hookrightarrow N_G(H)/H$. It follows
that $\Gamma_\n=\{1\}$. But since $\Gamma_\v$ is not connected
($\gamma\not\in \Gamma_\v^\circ$), we get $\Gamma_\n\neq\{1\}$.
\end{proof}

\subsection{Proof of Theorem \ref{Thm:1.1}}\label{SECTION_Proof1}
Assume the contrary, let $G,H_1,H_2$ be such that $H_1\equiv^G H_2$
but $H_1\not\sim_G H_2$.

In the proof we may assume that Theorem
\ref{Thm:1.1} is proved for all reductive groups $\underline{G}$
such that $\dim\underline{G}<\dim G$ and all subgroups
$\widetilde{H}_1,\widetilde{H}_2\subset G$ such that
$\widetilde{H}_1\equiv^G \widetilde{H}_2$ and
$\dim\widetilde{H}_1>\dim H_1$. Luna in \cite{Luna5}, Section 6,
proved that it is enough to prove Theorem \ref{Thm:1.1} only when
both $G/H_1,G/H_2$ have wonderful embeddings. His proof works for an
arbitrary reductive group $G$. In fact, he proved that, once Theorem 1 is proved for $G/H_1,G/H_2$
having wonderful embeddings (equivalently, $\X_{G,G/H_i}$ is spanned by $\Psi_{G,G/H_i}$,
see Proposition \ref{Prop:2.6}), homogeneous spaces $G/\underline{H}_1,G/\underline{H}_2$ with $\D_{G,G/\underline{H}_1}=\D_{G,G/\underline{H}_2},
\Psi_{G,G/\underline{H}_1}=\Psi_{G,G/\underline{H}_2}$ are isomorphic provided $\X_{G,G/\underline{H}_1}=
\X_{G,G/\underline{H}_2}$.

So it is enough to assume that both
$G/H_1,G/H_2$ admit wonderful embeddings. Fix some identification of
$\D_{G,G/H_1}$ and $\D_{G,G/H_2}$.

The following lemma allows to carry out an "induction step".

\begin{Lem}\label{Lem:3.1.1}
Let $\widetilde{H}\in \H_{H_1}$.  Then, possibly after conjugating
$H_2$ in $G$, we get $\widetilde{H}\in \H_{H_2}$,
$\D_{H_1}^{\widetilde{H}}=\D_{H_2}^{\widetilde{H}}$ and
$\iota_1(D)=\iota_2(D)$, where $\iota_i$ is the
 bijection $\D_{G,G/\widetilde{H}}\rightarrow \D_{G,G/H_i}\setminus
\D_{H_i}^{\widetilde{H}},i=1,2$, induced by taking the preimage.
\end{Lem}
\begin{proof}
Let $\widetilde{H}_2\in \H_{H_2}$  correspond to the colored
subspace $(\a_{H_1}^{\widetilde{H}},\D_{H_1}^{\widetilde{H}})$. By
Proposition \ref{Prop:2.8}, $\widetilde{H}\equiv^G \widetilde{H}_2$.
Since $\dim G/\widetilde{H}<\dim G/H_1$, it follows  that there is
$g\in G$ such that $\widetilde{H}=g\widetilde{H}_2 g^{-1}$. Replace
$H_2$ with $gH_2g^{-1}$. By Remark \ref{Rem:4.1}, the bijection
$\iota_2^{-1}\iota_1:\D_{G,G/\widetilde{H}}\rightarrow
\D_{G,G/\widetilde{H}}$ is induced by some element $g\in
N_G(\widetilde{H})$. To obtain $\iota_1(D)=\iota_2(D)$ it remains to
replace $H_2$ with $g^{-1}H_2g$. \end{proof}

The proof of Theorem 1 is in six steps. At first (step 1), we
reduce the proof to the case $N_G(H_i)=H_i,i=1,2$. Then (steps 2,3)
we show that it is enough to consider the situation when there is
$\widetilde{H}\in \underline{\H}_{H_1}\cap\underline{\H}_{H_2}$ and
$H_1,H_2,\widetilde{H}$ have a common maximal reductive subgroup,
say $S$. On steps 4,5 we show that one, in addition, can assume that
$\Rad_u(\h_1)\cong^S\Rad_u(\h_2)$ . Finally (step 6), we use
Proposition \ref{Cor:2.2}  to show that $H_1\sim_G H_2$.

{\it Step 1.} It follows from Theorem \ref{Thm:1.2} that
$N_G(H_1)\equiv_G N_G(H_2)$. By our assumptions, $G/H_i$ admits a
wonderful embedding. So $N_G(H_i)^\circ\subset H_i,i=1,2$ (this can
be deduced, for example, from Proposition \ref{Prop:2.6} and Lemma
\ref{Lem:2.1}). Assume that we have checked $N_G(H_1)\sim_G
N_G(H_2)$. Then Corollary \ref{Cor:2.3} implies $H_1\sim_G H_2$.
Proposition \ref{Prop:4.5} implies that at least one of subgroups
$H_1,H_2$ is contained in a proper parabolic subgroup of $G$.

{\it Step 2.} Let $Q^-$ be a proper parabolic subgroup of $G$ such
that $H_i\subset Q^-$ for some $i=1,2$. Let us show that
$\Rad_u(Q^-)\not\subset H_i$. Assume the contrary, let, say,
$\Rad_u(Q^-)\subset H_1\subset Q^-$. Set $M=Q^-/\Rad_u(Q^-)$. Thanks
to Lemma \ref{Lem:3.1.1}, we may assume that $H_2\subset Q^-$ and
$\D_{H_1}^{Q^-}=\D_{H_2}^{Q^-}$. Then, by Proposition
\ref{Prop:2.3}, $\Rad_u(Q^-)\subset H_2,H_1/\Rad_u(Q^-)\equiv^{M}
H_2/\Rad_u(Q^-)$. Since $\dim M<\dim G$, we see that $H_1\sim_{Q^-}
H_2$.

{\it Step 3.} Thanks to step 1, we may assume that there is $i\in
\{1,2\}$ and a proper parabolic subgroup $Q^-\in \H_{B^-}$ such that
$H_i\in \overline{\H}^{Q^-}$, and $S:=M\cap H_i$ is a maximal
reductive subgroup of $H_i$, where $M$ is the standard Levi subgroup
of $Q^-$. To be definite, assume $i=1$. By step 2 and assertion 2 of
Lemma \ref{Lem:5},  there exists $\widetilde{H}\in
\underline{\H}_{H_1}\cap\overline{\H}^{Q^-}$. We may replace $Q^-$
with a minimal parabolic subgroup of $G$ containing
$N_G(\widetilde{H})$. In particular, the subgroup $S\subset M$
satisfies the conditions of Proposition \ref{Prop:4.1}. Let
$\widetilde{H},H_1,H_2$ satisfy the conditions of Lemma
\ref{Lem:3.1.1} In this case, automatically,
$\D_{H_1}^{Q^-}=\D_{H_2}^{Q^-}$.

Clearly,  $\dim H_2/\Rad_u(H_2)\leqslant \dim S$. From symmetry
between $H_1$ and $H_2$ we get $\dim H_2/\Rad_u(H_2)=\dim S$ whence
$H_2\in \overline{\H}^{\widetilde{H}}$. Conjugating $H_2$ in
$\widetilde{H}$, we obtain that $S\subset H_2$.

{\it Step 4.} Let $Q^-,M,S$ be such  as on the previous step. To
simplify the notation, put $\n_i:=\Rad_u(\h_i),
\widetilde{\n}:=\Rad_u(\widetilde{\h}),\widetilde{\v}:=\Rad_u(\q^-)/\widetilde{\n},
\v_i:=\widetilde{\n}/\n_i,i=1,2,X:=G/\widetilde{H},
X':=M*_S\widetilde{\v}$.  Let us check that
$(\widetilde{\v}\oplus\v_1)\sim_{N_M(S)}(\widetilde{\v}\oplus\v_2)$.

Put $X_i':=M*_S(\widetilde{\v}\oplus \v_i),i=1,2$. By Lemma
\ref{Lem:2.12}, it is enough to check that $X_1'\cong^M X_2'$. From
Lemmas \ref{Lem:2.3}, \ref{Lem:2.4} it follows that $X_1'\equiv^M
X_2'$. Thanks to Proposition \ref{Prop:2.18} and the assumption that
Theorem \ref{Thm:1.1} holds for $M$, we get $X_1'\cong^M X_2'$.

{\it Step 5.}    Step 4 implies $\widetilde{H}\in
\underline{\H}_{H_2}$.

 Let us check that there is $g\in
N_G(\widetilde{H})\cap M$ such that  $\v_1\cong^S \v_2^g$.   Assume
that $\v_1\not\cong^S \v_2$.

Recall that $\Psi_{M,X'}=\Psi_{G,X}\cap\Span_\Q(\Delta(\m))$ and
$\A_{G,X}\subset \A_{M,X'}$ (Lemma \ref{Lem:2.3}). It follows from
Theorem \ref{Thm:1.2} and Lemma \ref{Lem:4.3} that
$\overline{\Psi}_{G,X}\cap\Span_\Q(\Delta(\m))=\overline{\Psi}_{M,X'}$.
So $\Lambda_{M,X'}$ is a direct summand of $\Lambda_{G,X}$.
Equivalently, $\A_{M,X'}/\A_{G,X}$ is connected. Therefore it is
enough to check that there is $g\in \A_{M,X'}$ such that
 $\v_1\cong^S \v_2^g$ (here we consider the natural homomorphism
 $\iota:\A_{M,X'}\rightarrow N_M(S)/S$). By step 4, there is $g\in
N_M(S)$ such that
\begin{equation}\label{eq:1}\widetilde{\v}\oplus\v_1\cong^S \widetilde{\v}^g\oplus\v_2^g.\end{equation}
 By assertion 2 of Proposition
\ref{Prop:4.1}, $g^2\in N_M(S)^\circ Z(M)^\circ$. Therefore
\begin{equation}\label{eq:2}
\widetilde{\v}^g\oplus \v_1^g\cong^S \widetilde{\v}\oplus \v_2.
\end{equation}

Adding (\ref{eq:1}) and (\ref{eq:2}), one easily gets
$\widetilde{\v}\cong^S \widetilde{\v}^g$ (in other words, $g$ is
lifted to $\A_{M,X'}$) and $\v_1\cong^S\v_2^g$.

Replacing $H_2$ with $gH_2g^{-1}$, we get $\v_1\cong^S \v_2$.

{\it Step 6.} Complete the proof. By Lemma \ref{Lem:5},
$[\widetilde{\n},\widetilde{\n}]\subset \n_1\cap\n_2$. Let $Y$ be
the subvariety of $\Gr_{\dim\n_i}(\widetilde{\n})$ consisting of all
$S$-stable subspaces $\n_0\subset \widetilde{\n}$ such that
$\n_1\cap\n_2\subset \n_0$ and $\n/\n_0\cong^S\n/\n_1$. Note that
$Y\cong \mathbb{P}^1$. Put $Y^0:=\{\n_0\in Y|
\n_\g(\s+\n_0)=\s+\n_0\}$. Since $\n_1,\n_2\in Y^0$, we get $Y^0\neq
\varnothing$. By Proposition \ref{Cor:2.2}, $\h_1\sim_G \h_2$.

\bigskip

{\Small Department of Mathematics, Massachusetts Institute of
Technology, Cambridge, MA 02139, USA.

E-mail address: ivanlosev@math.mit.edu}
\end{document}